\documentclass[10pt]{elsarticle}
\usepackage[cp1251]{inputenc}
\usepackage[english]{babel}
\usepackage{amsmath}
\usepackage{amssymb}
\usepackage{amsfonts}
\usepackage{amsthm}
\usepackage{mathtools}
\usepackage{dsfont}
\usepackage{rotating}
\usepackage{graphicx}
\usepackage{floatflt,epsfig}
\usepackage{lineno,hyperref}
\usepackage{enumerate}
\usepackage{colortbl}
\usepackage{array,tabularx,tabulary,booktabs}
\usepackage{longtable}
\usepackage{multirow}
\usepackage{wrapfig}
\usepackage{subcaption}
\usepackage{pdflscape}
\usepackage[table]{xcolor}
\usepackage[left=1in,right=1in,top=1in,bottom=1in]{geometry}
\newcolumntype{^}{>{\currentrowstyle}}

\journal{Arxiv}
\newtheorem{lem}{Lemma}[section]%

\newtheorem{prob}{Problem}%

\def\D{{\rm D}}

 \def\lg{\langle} \def\rg{\rangle}

\def\Aut{\hbox{\rm Aut\,}}  
  
  \def\mod{\hbox{\rm mod }}

  \def\GL{\hbox{\rm GL}}  \def\P\GL{\hbox{\rm P\GL}}
  \def\FF{{\hbox{\sf F\kern-.43emF}}}

\def\C{\hbox{\rm C}}

\def\ZZ{\mathbb{Z}} 

\def\nd{\mathrel{\bigm|\kern-.7em/}} 
 \def\f{\noindent}
\def\qed{\hfill $\Box$} \def\demo{\f {\bf Proof}\hskip10pt}

\begin{document}

\begin{frontmatter}
\title{Second largest maximal cliques in small Paley graphs of square order}

\author[01,02]{Huye Chen}
\ead{chenhy280@163.com}

\author[03]{Sergey Goryainov}
\ead{sergey.goryainov3@gmail.com}

\author[01]{Cong Hu}
\ead{18084314223@163.com}

\address[01] {Three Gorges Mathematical Research Center, China Three Gorges University, 8 University Avenue, Yichang 443002, Hubei Province, P.R. China}
\address[02] {School of Mathematics and Information Science, Guangxi University, Nanning, Guangxi, 530004, P. R. China}
\address[03] {School of Mathematical Sciences, Hebei International Joint Research Center for Mathematics and Interdisciplinary Science, Hebei Key Laboratory of Computational Mathematics and Applications, Hebei Workstation for Foreign Academicians,\\Hebei Normal University, Shijiazhuang  050024, P.R. China}


\begin{abstract}
There is a conjecture that the second largest maximal cliques in Paley graphs of square order $P(q^2)$ have size $\frac{q+\epsilon}{2}$, where $q \equiv \epsilon \pmod 4$, and split into two orbits under the full group of automorphisms whenever $q \ge 25$ (a symmetric description for these two orbits is known). However, some extra second largest maximal cliques (of this size) exist in $P(q^2)$ whenever $q \in \{9,11,13,17,19,23\}$. In this paper we analyse the algebraic and geometric structure of the extra cliques.
\end{abstract}

\begin{keyword}
 Paley graph, maximal cliques
\vspace{\baselineskip}
\MSC[2020] 05C25 \sep 11T30
\end{keyword}
\end{frontmatter}

\section{Introduction}
Let $q$ be an odd prime power such that $q\equiv 1 (\mod 4)$ , and let $\FF_{q}$ be the finite filed of order $q$ with primitive root $\delta$. The Paley graph of order $q$, denoted by $P(q)$, is an undirected graph defined on the elements of the finite field $\FF_{q}$ such that two vertices are adjacent if and only if the difference of these two vertices is a nonzero square in the multiplicative group $\FF_{q}^{\ast}$. The Paley graph $P(q)$ is known to be a self-complementary arc-transitive strongly regular graph with smallest eigenvalue $\frac{-1-\sqrt{q} }{2}$. According to the Delsarte-Hoffman bound \cite{P. Delsarte}, the size of a maximum independent set of a Paley graph $P(q)$ is at most $\sqrt{q}$. Since Paley graph $P(q)$ is self-complementary, the size of a maximum clique in $P(q)$ is also at most $\sqrt{q}$. 
 
 Through out this paper, set $\epsilon\in\{1,3\}$. Then $q^2\equiv 1 (\mod 4)$ for $q\equiv \epsilon (\mod 4)$. Now, we only study Paley graph of square order, which is denoted by $P(q^{2})$. For a Paley graph $P(q^{2})$, the Delsarte-Hoffman bound applied to $P(q^{2})$ gives that the size of a maximum independent set (a maximum clique) is $q$. It is obvious that the subfield $\FF_q$ of order $q$ forms a clique in Paley graph $P(q^{2})$. In 1984, Blokhuis \cite{A. Blokhuis} studied maximum cliques (of order $q$) in Paley graphs $P(q^{2})$ and proved that these cliques are isomorphic to the subfield $\FF_q$. In 1996, Baker et al. \cite{R.D. Baker} provided a  construction of maximal cliques of size $\frac{q+\epsilon}{2}$ for $q\equiv \epsilon (\mod 4)$. These cliques are maximal but not maximum and the cliques with this structure are not the only ones. They proposed a conjecture that there are no maximal cliques of size $s$ in $P(q^2)$, where $\frac{q+\epsilon}{2}< s<q$ for $q\equiv \epsilon (\mod 4)$. 
 In 2009, Kiermaier and Kurz \cite{Michael Kiermaier} considered integral point sets in affine planes over finite field and provided two maximal integral point sets and proved their maximality. In 2018, Goryainov et al. \cite{Sergey Goryainov} provided a construction of maximal cliques of size $\frac{q+\epsilon}{2}$ for $q\equiv \epsilon (\mod 4)$. By Magma \cite{Magma}, they found that for $25\le q\le83$, the graph $P(q^{2})$ contains exactly two non-equivalent (under the action of the automorphism group) maximal cliques of size $\frac{q+\epsilon}{2}$ for $q\equiv \epsilon (\mod 4)$, and these cliques are the second largest. In 2022, Goryainov et al. \cite{S. Goryainov} established a linear fractional correspondence between these two types of maximal cliques of size $\frac{q+\epsilon}{2}$ for $q\equiv \epsilon (\mod 4)$ in $P(q^{2})$. In 2024, Brouwer et al. \cite{S. Goryainov2} generalised the result by Baker et al. \cite{R.D. Baker} to the collinearity graphs of the Desargusian nets (in this paper we show that their construction gives maximal cliques of size $\frac{q+\epsilon}{2}$ in $P(q^2)$ for $q\in\{9,13,17\}$, which are constructed from two intersecting lines.

 In \cite{Sergey Goryainov}, the number of orbits on the maximal cliques of size $\frac{q+\epsilon}{2}$ in $P(q^2)$ was computed for $q\leq23$, and the results are summarised in the following table.
 \begin{table}[htp]
\centering
\begin{tabular}{cccccccccc}
\toprule  
$q$ & $3$ &$5$ &$7$ &$9$& $11$& $13$ & $17$ & $19$ & $23$ \\
\midrule  
Clique Size& $3$ &$3$ &$5$ &$5$& $7$& $7$ & $9$ & $11$ & $13$ \\
\midrule 
Number of Orbits& $1$ &$1$ &$1$ &$3$& $3$& $4$ & $9$ & $4$ & $4$ \\
\midrule 

\end{tabular}
\end{table}

\noindent
In this paper we analyse the algebraic and geometric structure of the extra cliques for $q \in \{9,11,13,17,19,23\}$ (except for the constructions from \cite{R.D. Baker} and \cite{Sergey Goryainov}). 

\section{Preliminaries}
In this section, we recall some basic facts about the affine planes and finite fields.
\subsection{Affine Plane $AG(2,q)$}
An affine plane is a point-line incidence structure. Let $AG(2,q)$ denote the affine plane, whose points are vectors of the $2$-dimensional vector $V(2,q)$ over $\FF_q$, and the lines are the additive shifts of $1$-dimensional subspaces of $V(2,q)$. It is well known that each line in $AG(2,q)$ contains $q$ points and there are $q+1$ lines through each point. There are $q^2$ points and $q(q+1)$ lines in $AG(2,q)$. 

\subsection{Finite fields}
Note that $\FF_{q^{2}}$ can be viewed as a quadratic extension over $\FF_{q}$. Set $\FF_q^*=\lg \delta \rg$ and let $\FF_{q^2}=\FF_q(\alpha)$ be the extension of $\FF_q$ by adding a root $\alpha$ of the irreducible polynomial $x^2-\delta$ over $\FF_q$. Then every element in $\FF_{q^{2}}$ can be uniquely written as $x+y\alpha$, for some $x,y \in\FF_q$. Moreover, $\FF_{q^2}$ can naturally be viewed as a $2$-dimensional vector spaces with basis $1$ and $\alpha$ over $\FF_{q}$. Define $\Phi$: $\FF_{q^{2}} \longrightarrow V(2,q)$ by $\Phi(x+y\alpha )=(x,y)$ for all $x$ and $y$ in $\FF_{q}$. It is easy to prove that $\Phi$ is a vector space isomorphism. Now we can identify the points of $AG(2,q)$ as the elements of $\FF_{q^{2}}$ and the line $l$ in $AG(2,q)$ can be viewed as $\left \{ x_1+y_1\alpha +c(x_2+y_2\alpha) \mid c\in \FF_q \right \}$, where $x_1+y_1\alpha\in\FF_{q^{2}}$ and $x_2+y_2\alpha\in \FF_{q^{2}}^{\ast}$ are fixed. We say a line $l$ of this form has a slope $(x_2+y_2\alpha)$ and note that the difference between any two points in $l$ is a scalar multiple of $(x_2+y_2\alpha)$. A line $l$ is called quadratic line (non-quadratic line) if the slope of this line is a square (non-square) element in $\FF_{q^{2}}^*$. For a quadratic line, the difference between any two distinct elements from this line is a square element in $\FF_{q^{2}}$. A Paley graph $P(q^{2})$ constructed on the points of $\FF_{q^{2}}$, and any two distinct points $a, b$ in $P(q^{2})$ are adjacent if and only if their difference is a square element in $\FF_{q^{2}}^{\ast}$. This implies that $a, b$ in $P(q^{2})$ are adjacent if and only if they are incident to a common line with square slope and the line through $a$ and $b$ is a quadratic line. Let $\beta$ be a primitive element in the finite field $\FF_{q^{2}}$ and $\delta$ be a primitive element in the subfield $\FF_{q}$. Note that $\FF_{q}^{\ast}=\left \langle \delta   \right \rangle =\left \langle \beta^{q+1}  \right \rangle$ and any element in $\FF_{q}^{\ast}$ is square in $\FF_{q^{2}}$. This further implies that the line $\FF_{q}$ is a clique of size $q$ in $P(q^{2})$.  

This must be throughout this paper, we set $S$ to be the set of square elements in $\FF_{q^2}^*$ and $S_0:=\{x+\alpha\mid x\in\FF_{q}\}\cap S$. Then we have that $S=\FF_q^{\ast}( S_0 \cup \{1\})$.
 \begin{lem}(\cite{S. Goryainov}, Proposition 1.)
 Through any point of $AG(2,q)$, there pass $q$ lines. Moreover, $\frac{q+1}{2}$ of these lines are quadratic and $\frac{q-1}{2}$ of these lines are non-quadratic. 
 \end{lem}

  \begin{lem}\label{aut}(\cite{G.A. Jones}, Theorem 9.1.)
  The automorphism group of Paley graph $P(q^{2})$ acts arc-transitively, and the equality
  $$A:=Aut(P(q^{2}))=\left \{ \gamma \mapsto a\gamma^{v}+b\mid a\in S,b\in \FF_{q^{2}},v\in Gal(\FF_{q^{2}})   \right \}$$ 
  holds, where $S$ is the set of square elements in $\FF_{q^{2}}^{\ast}$.

  The group $Aut(P(q^{2}))$ preserves the sets of quadratic and non-quadratic lines and acts on each of them transitively.
  \end{lem}

\begin{lem}\label{minusH}
Let $q$ be an odd prime power, $k$ be a divisor of $q-1$ and $H$ be the subgroup of order $k$ in $\mathbb{F}_q^*$. Then the following statements are hold.\\
{\rm (1)} If $k$ is even, then $-H = H$.\\
{\rm (2)} If $k$ is odd, then $-H \cap H = \emptyset$.
\end{lem}
\demo
{\rm (1)} If $k$ is even, for any $h\in H$, $(-h)^k=1$, then $-h\in H$. It follows that $-H = H$.\\
{\rm (2)} On the contrary, assume that there exists an element $h \in -H \cap H$. Then $h^k=1$ and $-h\in H$, this implies $(-h)^k=1$, a contradiction. 
\qed

\vskip 5mm
Followed by the above notations, the following lemma is obvious.
\begin{lem} (\cite{S. Goryainov}, Proposition 2.)
 {\rm (1)} The element $-1$ is a square in $\FF_{q}^{\ast}$ if and only if $q\equiv1 (\mod 4)$.\\
  {\rm (2)} The element $\alpha$ is a square in $\FF_{q^2}^{\ast}$ if and only if $q\equiv3 (\mod 4)$.
\end{lem} 

\section{Constructions}
 In this section, we propose the constructions of maximal cliques in Paley graph $P(q^{2})$ of size $\frac{q+\epsilon}{2}$ for $q\equiv \epsilon (\mod 4)$, where $9\le q\le 23$. Moreover, we give some numerical information of these cliques by computing the stabilizer of each clique under the action of the automorphism group $Aut(P(q^2))$. In this section, we set $\D_n$ denote the dihedral group of order $n$, and set $\ZZ_n$ denote the cyclic group of order $n$. Throughout this section, by $Aut(P(q^2))_C$ denote the stabilizer of the automorphism group $Aut(P(q^2))$ acting on the set of cliques, where $C$ is a clique of the Paley graph $P(q^2)$.

\subsection{Maximal cliques in Paley graph $P(9^2)$}

\begin{figure}[htp]
    \centering
    \includegraphics[height=4cm, width=8cm]{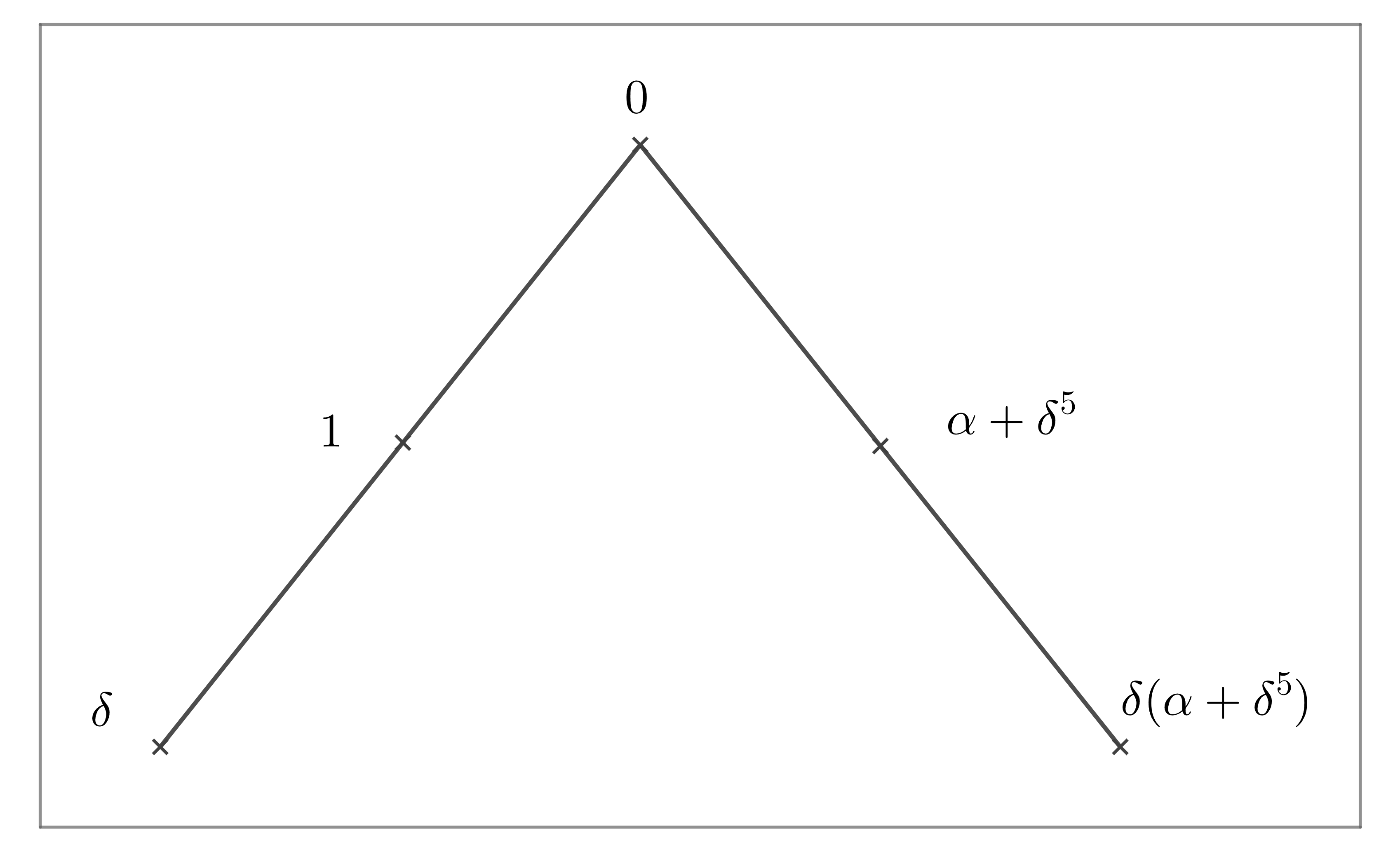}
    \quad 
    \caption{\label{$C_9$} $C_9$}
\end{figure}

We can choose a primitive element $\delta\in \FF_9$ and $\alpha\in\FF_{9^2}$ be a root of   the irreducible polynomial $x^2-\delta$ over $\FF_q$,
such that $S_0=\{1, \alpha +\delta,\alpha +\delta^{2},\alpha +\delta^{5},\alpha +\delta^{6}\}$.
Let $H= \{ 1,\delta \}$  be a subset in $\FF_{9}^{\ast}$. 
Set $C_9:=\{0\}\cup H\cup (\alpha+\delta^5)H$ be a subset of finite field $\FF_{9^2}$. By Magma, we have that $C_9$ is a maximal clique in $P(9^2)$ with size $\frac{q+1}{2}$ for $q=9$.
The structure of the clique $C_9$ is presented in Fig \ref{$C_9$}.

\begin{lem}
Set $\mathcal{C}_9$ be the orbit of $Aut(P(9^2))$ acting on the cliques with $C_9\in \mathcal{C}_9$. Then $Aut(P(9^2))_{C_9}=\lg \phi' \rg\cong \ZZ_2$, where $\phi'(\gamma)=(\alpha+\delta^5)\gamma^9$ for any $\gamma\in \FF_{9^2}$. Moreover, $|\mathcal{C}_9|=6480$.
\end{lem}
\demo
Note that for any $\phi \in \Aut(P(9^2))$,
$\phi(\gamma)=a\gamma^{v}+b$, where $ a\in S$, $b\in \FF_{9^{2}}$, $v\in Gal(\FF_{9^{2}})$.  If $\phi\in Aut(P(9^{2}))_{C_{9}}$, then  $\phi(0)=0$ and $b=0$. Set $\mathcal{H}:=\{H,(\alpha+\delta^5)H\}$ be the subset of two intersection lines, which are presented in Fig \ref{$C_9$}. Then $\phi(\mathcal{H})=\mathcal{H}$. It follows that either $\phi(H)=H$, $\phi((\alpha +\delta ^{5})H)=(\alpha +\delta ^{5})H$ or $\phi(H)=(\alpha +\delta ^{5})H$, $\phi((\alpha +\delta ^{5})H)=H$.

If $\phi(H)=aH^v=H$ and $\phi((\alpha +\delta ^{5})H)=a(\alpha +\delta ^{5})^vH^v=(\alpha +\delta ^{5})H$, then we have that $\{a,a\delta^v\}=\{1,\delta\}$ and $\{a(\alpha+\delta^5)^v,a\delta^v(\alpha+\delta^5)^v\}=\{\alpha+\delta^5,\delta(\alpha+\delta^5)\}$. It follows that $\phi(\gamma)=\gamma$ for any $\gamma\in \FF_{9^2}$. Now, $\phi$ is the identity automorphism.

If $\phi'(H)=(\alpha +\delta ^{5})H$ and  $\phi'((\alpha +\delta ^{5})H)=H$, then by the same arguments as in the first case we have that $a=\alpha +\delta ^{5}$ and $|v|=2$. Then we have that $\phi'(\gamma)=(\alpha +\delta ^{5})\gamma^{9}$ for any $\gamma\in \FF_{9^2}$. It follows that $\phi'^2=1$.

Now, we have that $Aut(P(9^{2}))_{C_{9}}=\{\phi,\phi'\}=\lg \phi' \rg \cong \ZZ_2$.
The action of $\phi'$ on the points in clique $C_{9}$ is presented in the following table.
 \begin{table}[htp]
\centering
\begin{tabular}{cccccc}
\toprule  
$\gamma$ & 0 &$1$ &$\delta$ &$\alpha +\delta ^{5}$& $\delta(\alpha +\delta ^{5})$\\
	\midrule  
$\phi'(\gamma)$ & 0 &$\alpha +\delta ^{5}$ &$\delta(\alpha +\delta ^{5})$ &1& $\delta$\\
\midrule 
\end{tabular}
\end{table}

Note that $\left | Aut(P(9^{2})) \right |=\frac{q^2-1}{2}\times q^{2}\times 4=12960$, so we have that $|\mathcal{C}_{9}|=|Aut(P(9^2)):Aut(P(9^2))_{C_9}|=6480$. 
\qed

\vskip 5mm

Set $C_{s,\gamma}^{\eta}:=\{sx^{\eta}+\gamma| x\in C_9\cap \FF_{9^2}\}$ for $\eta\in\{1,3\}$, where $s\in S$ and $\gamma\in\FF_{9^2}$ . 

\begin{lem}
Set $\mathcal{C}_9$ be the orbit of $Aut(P(9^2))$ acting on the cliques with $C_9\in \mathcal{C}_9$. Then $\mathcal{C}_9=\{C_{s,\gamma}\mid s\in S, \gamma\in\FF_{9^2}\}\cup \{C_{s,\gamma}^3\mid s\in S, \gamma\in\FF_{9^2}\}$.
\end{lem}
\demo
It is obvious that $\{C_{s,\gamma}\mid s\in S, \gamma\in\FF_{9^2}\}\cup \{C_{s,\gamma}^3\mid s\in S, \gamma\in\FF_{9^2}\}\subset \mathcal{C}_9$. Now we will prove that $C_{s,\gamma}^{\eta}=C_{s',\gamma'}^{\eta'}$ if and only if $s=s'$, $\gamma=\gamma'$ and $\eta=\eta'$, where $s,s'\in S$, $\gamma,\gamma'\in\FF_{9^2}$ and $\eta,\eta'\in\{1,3\}$  .

Note that the intersection point of two lines in the clique $C_{s,\gamma}^{\eta}$  is $\gamma$. If $C_{s,\gamma}^{\eta}=C_{s',\gamma'}^{\eta'}$, then $\gamma=\gamma'$ and $\{s,s\delta^{\eta}\}\in\{\{s',s'\delta^{\eta'}\},\{s'(\alpha+\delta^5)^{\eta'},s'\delta^{\eta'}(\alpha+\delta^5)^{\eta'}\}\}$.

If $\{s,s\delta^{\eta}\}=\{s',s'\delta^{\eta'}\}$, then $s=s'$ and $\eta=\eta'$.

If $\{s,s\delta^{\eta}\}=\{s'(\alpha+\delta^5)^{\eta'},s'\delta^{\eta'}(\alpha+\delta^5)^{\eta'}\}$, then 
either $s=s'(\alpha+\delta^5)^{\eta'}$, 
$s\delta^{\eta}=s'\delta^{\eta'}(\alpha+\delta^5)^{\eta'}$ or
$s=s'\delta^{\eta'}(\alpha+\delta^5)^{\eta'}$, 
$s\delta^{\eta}=s'(\alpha+\delta^5)^{\eta'}$.

Case 1:
If $s=s'(\alpha+\delta^5)^{\eta'}$ and
$s\delta^{\eta}=s'\delta^{\eta'}(\alpha+\delta^5)^{\eta'}$, then $s'(\alpha+\delta^5)^{\eta'}\delta^{\eta}=s'\delta^{\eta'}(\alpha+\delta^5)^{\eta'}$ and 
$\eta=\eta'$.
Now $\{s(\alpha+\delta^5)^{\eta},s\delta^{\eta}(\alpha+\delta^5)^{\eta}\}=\{s',s'\delta^{\eta'}\}$ with $\eta=\eta'$, then 
either $s(\alpha+\delta^5)^{\eta}=s'$ or $s(\alpha+\delta^5)^{\eta}=s'\delta^{\eta'}$.
If $s(\alpha+\delta^5)^{\eta}=s'$, then 
$s=s(\alpha+\delta^5)^{\eta}(\alpha+\delta^5)^{\eta}$, a contradiction. If $s(\alpha+\delta^5)^{\eta}=s'\delta^{\eta}$, then we have  that $(\alpha+\delta^5)^{\eta^2}=\delta^{\eta}$ followed from $s=s'(\alpha+\delta^5)^{\eta'}$, a contradiction.

Case 2:
If $s=s'\delta^{\eta'}(\alpha+\delta^5)^{\eta'}$, 
$s\delta^{\eta}=s'(\alpha+\delta^5)^{\eta'}$, then we have that $\delta^{\eta+\eta'}=1$, a contradiction.

Followed from the above arguments, we have that $\mathcal{C}_9=\{C_{s,\gamma}\mid s\in S, \gamma\in\FF_{9^2}\}\cup \{C_{s,\gamma}^3\mid s\in S, \gamma\in\FF_{9^2}\}$, because $|\{C_{s,\gamma}\mid s\in S, \gamma\in\FF_{9^2}\}\cup \{C_{s,\gamma}^3\mid s\in S, \gamma\in\FF_{9^2}\}|=\frac{q^2-1}{2}\times q^{2}\times 2=6480$.
\qed

\subsection{Maximal cliques in Paley graph $P(13^2)$}
We can choose a primitive element $\delta\in \FF_{13}$ and $\alpha\in\FF_{{13}^2}$ be a root of  the irreducible polynomial $x^2-\delta$ over $\FF_q$,
such that $\FF_{13}^*=\lg\delta\rg=\lg 6 \rg$ and $S_0=\{1,\alpha+8, \alpha+5, \alpha+7, \alpha+6, \alpha+1, \alpha+12\}$.
\subsubsection{$C_{13}^{A}$-construction}

Let $H=\left \{ \delta^{4} ,\delta ^{8},\delta ^{12}  \right \}=\left \{ 9,3,1  \right \} $ be a subgroup of $\FF_{13}^*$ with order $3$. 
Then 
$\FF_{13}^{\ast}=H \cup \delta H \cup \delta ^{2} H \cup \delta ^{3} H=\left \{ 9,3,1  \right \} \cup  \left \{10,12,4  \right \}\cup  \left \{ 8,7,11 \right \}\cup \left \{ 6,2,5 \right \}$.

\begin{figure}[htp]
    \centering
    \includegraphics[height=4cm, width=8cm]{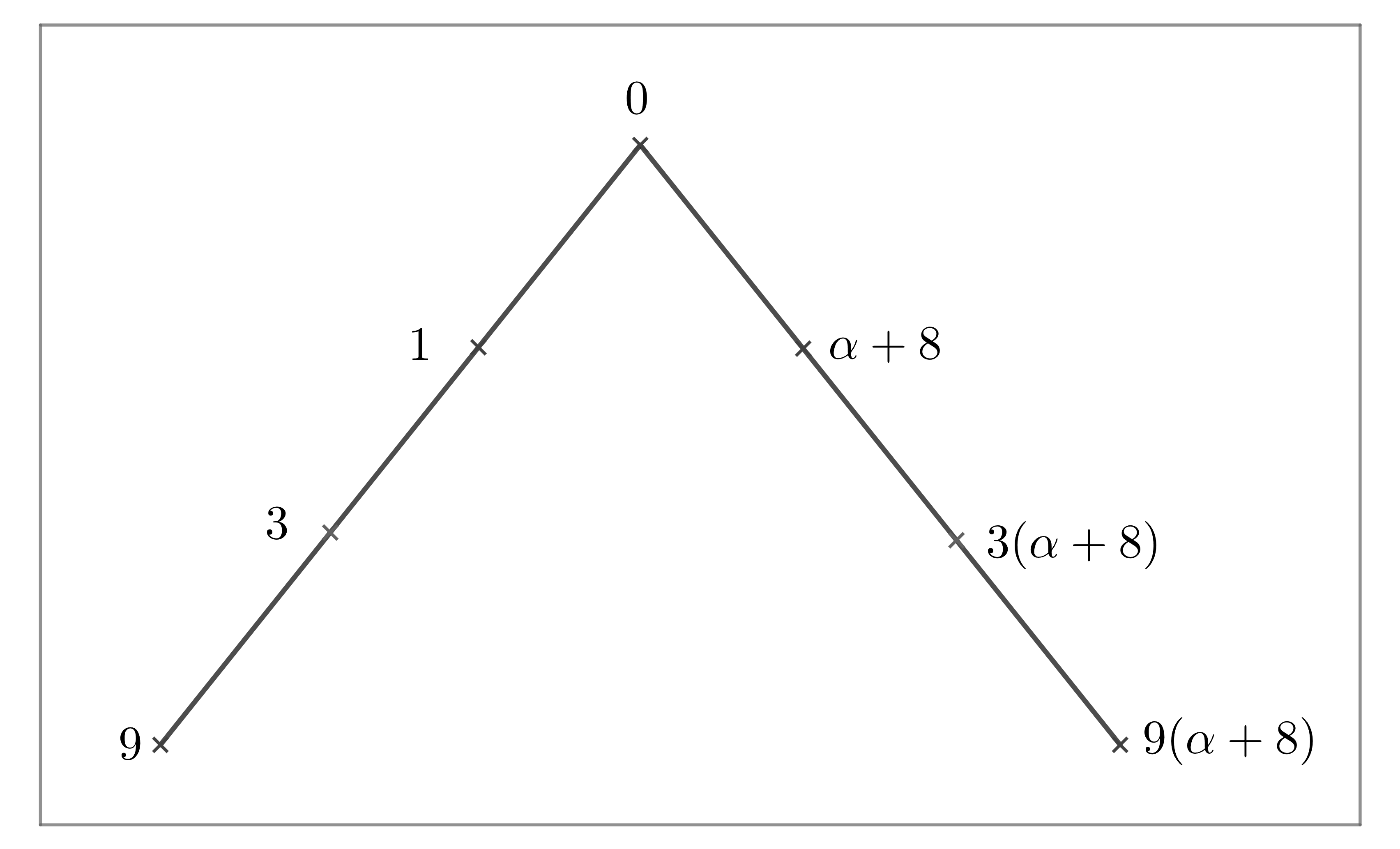}
    \quad
    \caption{\label{q=13-1}$C_{13}^{A}$}
    \end{figure}

Set $C_{13}^{A}:=\{0\}\cup H\cup (\alpha+8)H$ be a subset of finite field $\FF_{13^2}$. By Magma, we have that $C_{13}^{A}$ is a maximal clique in $P(13^2)$ with size $\frac{q+1}{2}$ for $q=13$.
The structure of the clique $C_{13}^{A}$ is presented in Fig \ref{q=13-1}.

\begin{lem}
Set $\mathcal{C}_{13}^{A}$ be the orbit of $Aut(P(13^2))$ acting on the cliques with $C_{13}^{A}\in \mathcal{C}_{13}^{A}$. Then $Aut(P(13^2))_{C_{13}^{A}}=\lg \sigma,\tau \rg\cong \D_6$, where $\sigma(\gamma)=3\gamma$ and $\tau(\gamma)=(\alpha+8)\gamma^{13}$ for any $\gamma\in \FF_{13^2}$.  Moreover, $|\mathcal{C}_{13}^{A}|=4732$.
\end{lem}
\demo
Note that for any $\phi \in \Aut(P(13^2))$,
$\phi(\gamma)=a\gamma^{v}+b$, where $ a\in S$, $b\in \FF_{13^{2}}$, $v\in Gal(\FF_{13^{2}})$. If $\phi\in Aut(P(13^{2}))_{C_{13}^{A}}$, then  $\phi(0)=0$ and $b=0$. Set $\mathcal{H}:=\{H,(\alpha+8)H\}$ be the subset of two lines, which are presented in Fig \ref{q=13-1}. Then $\phi(\mathcal{H})=\mathcal{H}$. It follows that either $\phi(H)=H$, $\phi((\alpha +8)H)=(\alpha +8)H$ or $\phi(H)=(\alpha +8)H$, $\phi((\alpha +8)H)=H$.

If $\phi(H)=H$, $\phi((\alpha +8)H)=(\alpha +8)H$, then we have that $aH^{v}=aH=H$ and $a((\alpha+8)H)^{v}=a(\alpha+8)^{v}H=(\alpha+8)H$. It follows that $a\in H$ and $v=1$. So $\phi(\gamma)=a\gamma$ with $a\in H$.

If $\phi(H)=(\alpha +8)H$, $\phi((\alpha +8)H)=H$, then we have that $aH^{v}=aH=(\alpha+8)H$ and $a((\alpha+8)H)^{v}=a(\alpha+8)^{v}H=H$, then $a\in (\alpha+8)H$ and $|v|=2$. So $\phi(\gamma)=a\gamma^{13}$ with $a\in (\alpha+8)H$. 

Let $\phi_1(\gamma)=\gamma$, $\phi_2(\gamma)=3\gamma$, $\phi_3(\gamma)=9\gamma$, $\phi_1'(\gamma)=(\alpha+8)\gamma^{13}$, $\phi_2'(\gamma)=3(\alpha+8)\gamma^{13}$, $\phi_3'(\gamma)=9(\alpha+8)\gamma^{13}$ for any $\gamma\in \FF_{13^2}$. Set $\sigma=\phi_2$ and $\tau=\phi_1'$. Then $\sigma^{3}=1$, $\tau^{2}=1$, $\tau^{-1}\sigma\tau=\sigma^{-1}$. Then $\lg \sigma,\tau \rg=\{1, \sigma, \sigma^2, \tau, \tau\sigma, \tau\sigma^2\}=\{\phi_1, \phi_2, \phi_3,\phi_1',\phi_2',\phi_3'\}$. 
The action of $\sigma$ and $\tau$ on the points in clique $C_{13}^{A}$ is presented by the following table. 
\vskip 1mm

 \begin{table}[htp]
\centering
\begin{tabular}{cccccccc}
\toprule  
$\gamma$ & 0 &$1$ &$3$ &$9$& $\alpha+8$& $3(\alpha+8)$& $9(\alpha+8)$\\
	\midrule  
$\sigma(\gamma)$ & 0 &$3$ &$9$ &1& $3(\alpha+8)$& $9(\alpha+8)$& $\alpha+8$\\
$\tau(\gamma)$ & 0 & $\alpha +8$ & $3(\alpha+8)$ &$9(\alpha+8)$ & $1$ &$3$ &$9$ \\
\midrule 
\end{tabular}
\end{table}

Followed by the above arguments, we know that $\sigma(C_{13}^A)=C_{13}^A$ and  $\tau(C_{13}^A)=C_{13}^A$. Then $Aut(P(13^2))_{C_{13}^{A}}=\lg \sigma,\tau\rg\cong \D_{6}$, where $\sigma(\gamma)=3\gamma$, $\tau(\gamma)=(\alpha+8)\gamma^{13}$ for any $\gamma\in \FF_{13^2}$.

Note that $\left | Aut(P(13^{2})) \right |=\frac{q^2-1}{2}\times q^{2}\times 2=28392$, so we have that $|\mathcal{C}_{13}^{A}|=|Aut(P(13^2)):Aut(P(13^2))_{C_{13}^{A}}|=4732$. 
\qed

\vskip 5mm
Set $C_{s,\gamma,i}:=\left \{ s\delta ^{i}x +\gamma   \mid  x\in C_{13}^{A}\cap \FF_{13^2}  \right \} $, where $s\in S_{0}$, $i\in \{\ 0,1,2,3 \ \}$ and  $\gamma \in \FF_{13^{2}}$. 

\begin{lem}
Set $\mathcal{C}_{13}^{A}$ be the orbit of $Aut(P(13^2))$ acting on the cliques with $C_{13}^{A}\in \mathcal{C}_{13}$. Then  $\mathcal{C}_{13}^{A}=\{C_{s,\gamma,i}\mid s\in S_{0},i\in \{0,1,2,3\}, \gamma\in\FF_{13^2}\}$.
\end{lem}
\demo
 It is obvious that $\{C_{s,\gamma,i}\mid s\in S_{0},i\in \{0,1,2,3\}, \gamma\in\FF_{13^2}\} \subset \mathcal{C}_{13}^{A}$. Now we will prove that $\C_{s,\gamma,i}=C_{s',\gamma',i'}$ if and only if $s=s'$, $\gamma=\gamma'$ and $i=i'$, where $s,s'\in S_0$, $i,i'\in\{0,1,2,3\}$ and $\gamma,\gamma'\in\FF_{13^2}$.

Note that the intersection point of two lines in the clique $C_{s,\gamma,i}$ is $\gamma$. If $C_{s,\gamma,i}=C_{s',\gamma',i'}$, then $\gamma=\gamma'$ and $\{s\delta^{i}H,s\delta^{i}(\alpha+8)H\}=\{s'\delta^{i'}H,s'\delta^{i'}(\alpha+8)H\}$. 

If $s\delta^{i}H=s'\delta^{i'}H$ and $s\delta^{i}(\alpha+8)H=s'\delta^{i'}(\alpha+8)H$, then $s\delta^i\in s'\delta^{i'}H$. If $s\delta^i= s'\delta^{i'}$, then $s=s'$ and $i=i'$; if $s\delta^i= \delta^4 s'\delta^{i'}$, then $\delta^i=\delta^{4+i'}$, a contradiction; if $s\delta^i= \delta^8 s'\delta^{i'}$, then $\delta^i=\delta^{8+i'}$, a contradiction.

If $s\delta^{i}H=s'\delta^{i'}(\alpha+8)H$ and $s'\delta^{i'}H=s\delta^{i}(\alpha+8)H$, then $(\alpha+8)^2 \in H$, a contradiction. 

Followed from the above arguments, we have that  $\mathcal{C}_{13}^{A}=\{C_{s,\gamma,i}\mid s\in S_{0}, i\in \{0,1,2,3\}, \gamma\in\FF_{13^2}\}$, because $|\{C_{s,\gamma,i}\mid s\in S_{0},i\in \{0,1,2,3\}, \gamma\in\FF_{13^2}\}|=\frac{q+1}{2}\times q^{2}\times 4=4732$.
\qed

\subsubsection{$C_{13}^{B}$-construction}
\begin{figure}[htp]
    \centering
    \includegraphics[height=4cm, width=8cm]{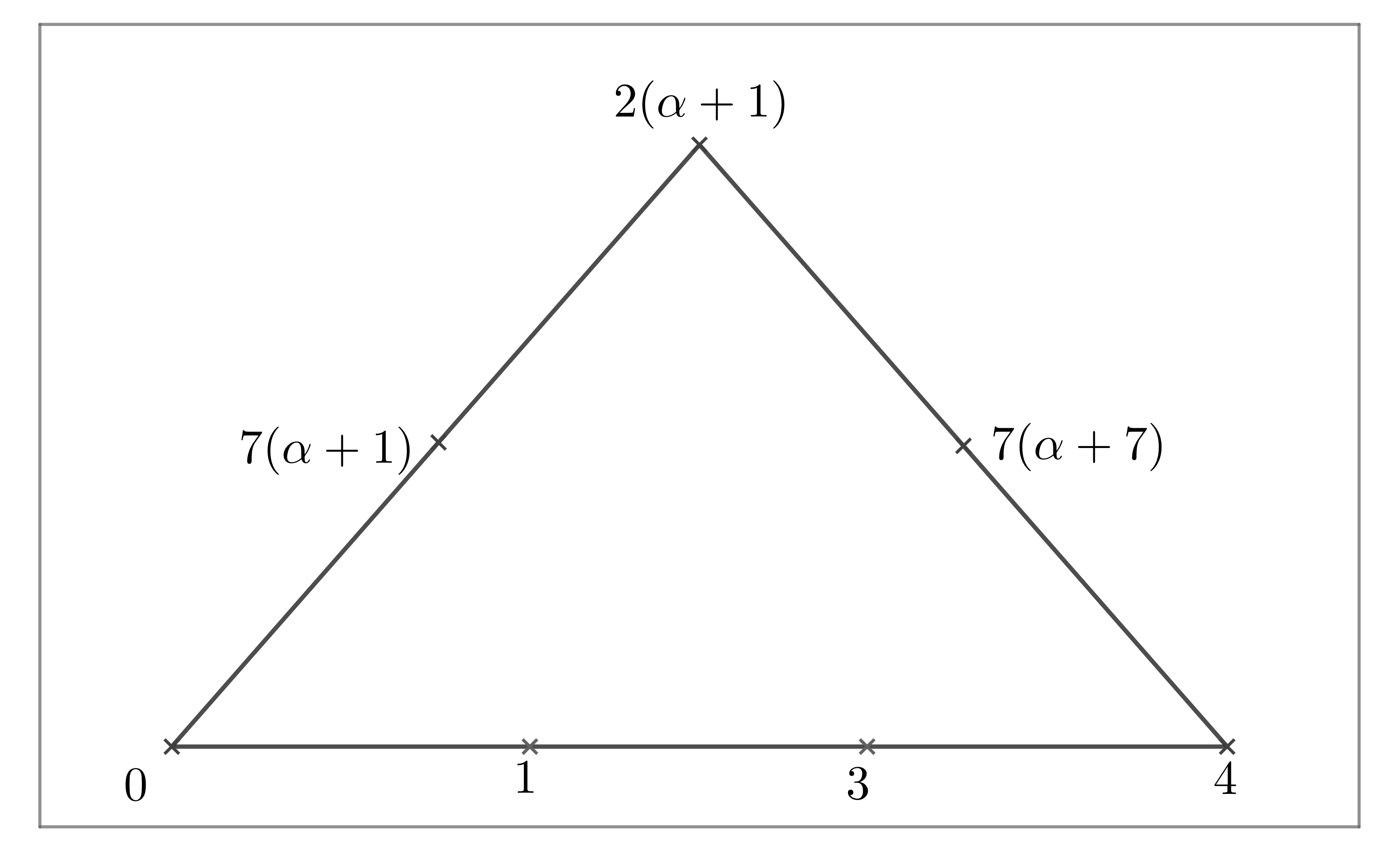}
    \quad
     \caption{\label{q=13-2}$C_{13}^{B}$}
    \end{figure}

Let $H= \{ 1,3,4 \}$  be a subset in $\FF_{13}^{\ast}$. Set $C_{13}^{B}:=\{0\}\cup H\cup\{7(\alpha+1), 2(\alpha+1), 7(\alpha+7)\}$ be a subset of finite field $\FF_{13^2}$. By Magma, we have that $C_{13}^{B}$ is a maximal clique in $P(13^2)$ with size $\frac{q+1}{2}$ for $q=13$.
The structure of the clique $C_{13}^{B}$ is presented in Fig \ref{q=13-2}.

\begin{lem}
Set $\mathcal{C}_{13}^{B}$ be the orbit of $Aut(P(13^2))$ acting on the cliques with $C_{13}^{B}\in \mathcal{C}_{13}^{B}$. Then $Aut(P(13^2))_{C_{13}^{B}}=\lg \phi' \rg\cong \ZZ_2$, where $\phi'(\gamma)=12\gamma^{13}+4$ for any $\gamma\in \FF_{13^2}$.  Moreover,  $|\mathcal{C}_{13}^{B}|=14196$.
\end{lem}

\demo
Note that for any $\phi \in \Aut(P(13^2))$,
$\phi(\gamma)=a\gamma^{v}+b$, where $ a\in S$, $b\in \FF_{13^{2}}$, $v\in Gal(\FF_{13^{2}})$.  If $\phi\in Aut(P(13^{2}))_{C_{13}^B}$, then  $\phi(2(\alpha+1))=2(\alpha+1)$. Set $\mathcal{P}:=\{0,4\}$ be the two special intersecting points in clique $C_{13}^{B}$, then $\phi(\mathcal{P})=\mathcal{P}$. It follows that either $\phi(0)=0$, $\phi(4)=4$ or $\phi(0)=4$, $\phi(4)=0$.

If $\phi(0)=0$ and $\phi(4)=4$, then we have that $a=1$, $b=0$. Then $\phi(\gamma)=\gamma^v$, $\phi(7(\alpha+1))=7(\alpha+1)$ and $\phi(7(\alpha+7))=7(\alpha+7)$. It follows that $\phi(\gamma)=\gamma$ for any $\gamma\in \FF_{13^2}$.

If $\phi'(0)=4$ and $\phi'(4)=0$. Then $a=12$, $b=4$ and $\phi'(\gamma)=12\gamma^v+4$. Moreover, we have that $\phi'(7(\alpha+1))=7(\alpha+7)$ and $\phi'(7(\alpha+7))=7(\alpha+1)$. Then we have that $\phi'(\gamma)=12\gamma^{13}+4$ for any $\gamma\in \FF_{13^2}$ and $\phi'^2=1$.

The action of $\phi'$ on the point in clique $C_{13}^{B}$ is presented by the following table.

 \begin{table}[htp]
\centering
\begin{tabular}{cccccccc}
\toprule  
$\gamma$ & 0 &$1$ &$3$ &$4$& $7(\alpha+1)$& $2(\alpha+1)$& $7(\alpha+7)$\\
	\midrule  
$\phi'(\gamma)$ & 4 & $3$ & $1$ & 0 & $7(\alpha+7)$& $2(\alpha+1)$& $7(\alpha+1)$\\
\midrule 
\end{tabular}
\end{table}

Followed by the above arguments, we get that $\phi'(C_{13}^B)=C_{13}^B$ and $Aut(P(13^2))_{C_{13}^{B}}=\lg \phi' \rg\cong \ZZ_2$, where $\phi'(\gamma)=12\gamma^{13}+4$ for any $\gamma\in \FF_{13^2}$.
Note that $\left | Aut(P(13^{2})) \right |=\frac{q^2-1}{2}\times q^{2}\times 2=28392$, so we have that $|\mathcal{C}_{13}^{B}|=|Aut(P(13^2)):Aut(P(13^2))_{C_{13}^{B}}|=14196$. 
\qed

\vskip 5mm

Set $C_{s,\gamma}:=\{sx+\gamma| x\in C_{13}^{B}\cap \FF_{13^2}\}$ where $s\in S$ and $\gamma\in\FF_{13^2}$.

\begin{lem}

Set $\mathcal{C}_{13}^{B}$ be the orbit of $Aut(P(13^2))$ acting on the cliques with $C_{13}^{B}\in \mathcal{C}_{13}^{B}$. Then $\mathcal{C}_{13}^{B}=\{C_{s,\gamma}\mid s\in S, \gamma\in\FF_{13^2}\}$.

\end{lem}
\demo
It is obvious that $\{C_{s,\gamma}\mid s\in S, \gamma\in\FF_{13^2}\}\subset \mathcal{C}_{13}$. Now we will prove that $C_{s,\gamma}=C_{s',\gamma'}$ if and only if $s=s'$ and $\gamma=\gamma'$  where $s,s'\in S$, $\gamma,\gamma'\in\FF_{13^2}$ .

If $C_{s,\gamma}=C_{s',\gamma'}$, then $2(\alpha+1)s+\gamma=2(\alpha+1)s'+\gamma'$ and the two special intersecting points $\{\gamma, 4s+\gamma\}=\{\gamma',4s'+\gamma\}$. If $\gamma=4s'+\gamma'$ and $\gamma'=4s+\gamma$, then $(\alpha+1)=1$, a contradiction. Now, we have $s=s'$ and $\gamma =\gamma '$.
And then $\mathcal{C}_{13}^{B}=\{C_{s,\gamma}\mid s\in S, \gamma\in\FF_{13^2}\}$, because $|\{C_{s,\gamma}\mid s\in S, \gamma\in\FF_{13^2}\}|=\frac{q^2-1}{2}\times q^{2} =14196$.
\qed

\subsection{Maximal cliques in Paley graph $P(17^2)$}
We can choose a primitive element $\delta\in \FF_{17}$ and $\alpha\in\FF_{{17}^2}$ be a root of  the irreducible polynomial $x^2-\delta$ over $\FF_q$,
such that $\FF_{17}^*=\lg\delta\rg=\lg 10 \rg$ and $S_0=\{1,\alpha+7, \alpha+10, \alpha+14, \alpha+3, \alpha+6, \alpha+11, \alpha+8, \alpha+9\}$.
\subsubsection{$C_{17}^{A}$-construction}
Let $H=\left \{ \delta^{4} ,\delta ^{8},\delta ^{12} , \delta ^{16} \right \}=\left \{ 4,16,13,1 \right \} $ be a subgroup in $\FF_{17}^{\ast}$. Then
$\FF_{q}^{\ast}=H \cup \delta H \cup \delta ^{2} H \cup \delta ^{3} H=\left \{ 4,16,13,1 \right \}\cup  \left \{6,7,11,10  \right \}\cup  \left \{ 9,2,8,5\right \}\cup\left \{ 5,3,12,14  \right \}$.

\begin{figure}[htp]
    \centering
    \includegraphics[height=4cm, width=8cm]{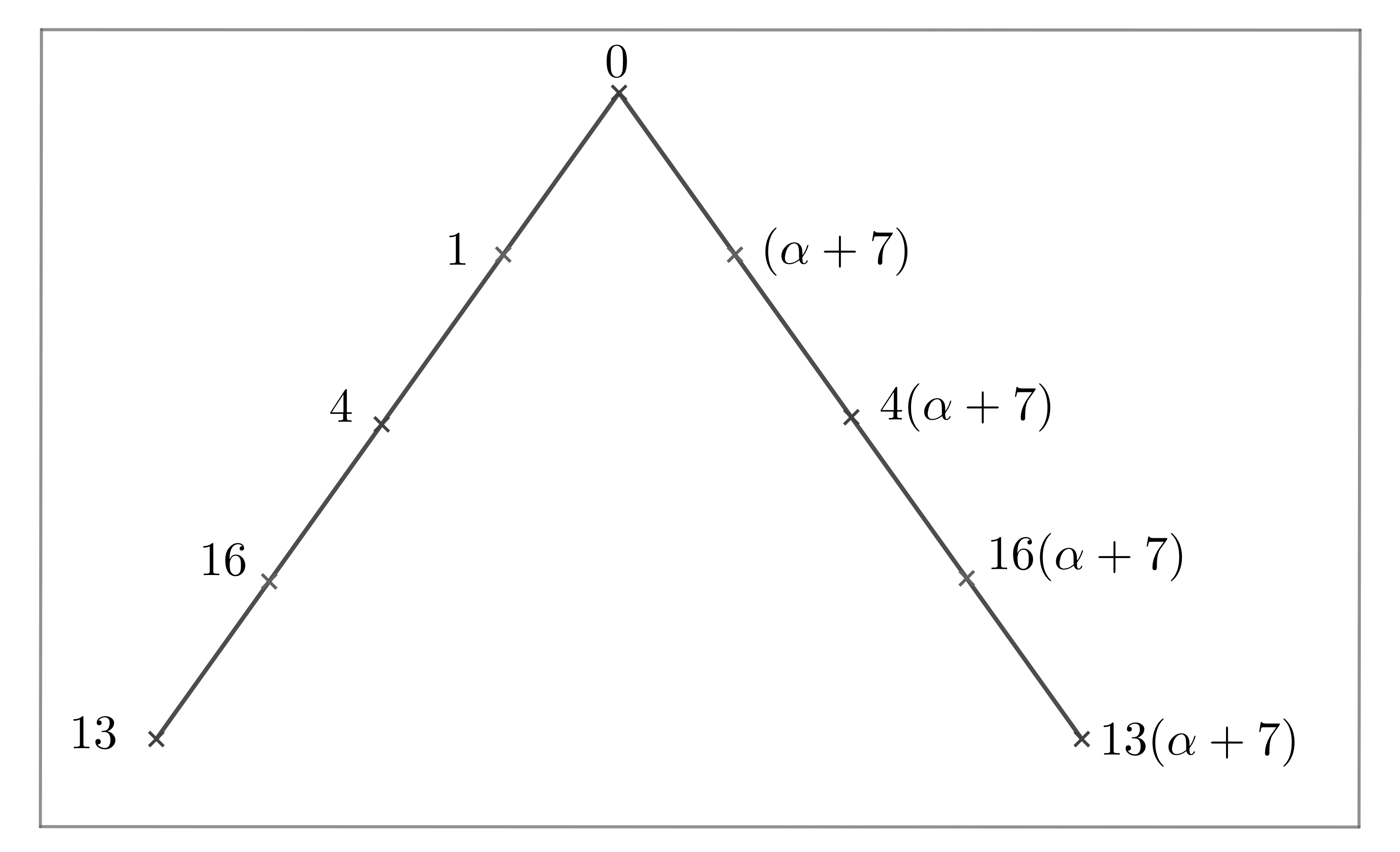}
    \quad
     \caption{\label{q=17-1}$C_{17}^{A}$}
    \end{figure}

Set $C_{17}^{A}:=\{0\}\cup H\cup (\alpha+7)H$ be a subset of finite field $\FF_{17^2}$. By Magma, we have that $C_{17}^{A}$ is a maximal clique in $P(17^2)$ with size $\frac{q+1}{2}$ for $q=17$.
The structure of the clique $C_{17}^{A}$ is presented in Fig \ref{q=17-1}.

\begin{lem}
Set $\mathcal{C}_{17}^{A}$ be the orbit of $Aut(P(17^2))$ acting on the cliques with $C_{17}^{A}\in \mathcal{C}_{17}^{A}$. Then $Aut(P(17^2))_{C_{17}^{A}}=\lg \sigma \rg\cong \ZZ_{8}$, where $\sigma(\gamma)=(\alpha+7)\gamma^{17}$ for any $\gamma\in \FF_{17^2}$. Moreover, $|\mathcal{C}_{17}^{A}|=10404$.
\end{lem}
\demo
Note that for any $\phi \in \Aut(P(17^2))$,
$\phi(\gamma)=a\gamma^{v}+b$, where $ a\in S$, $b\in \FF_{17^{2}}$, $v\in Gal(\FF_{17^{2}})$. If $\phi\in Aut(P(17^{2}))_{C_{17}^{A}}$, then  $\phi(0)=0$ and $b=0$. Set $\mathcal{H}:=\{H,(\alpha+7)H\}$ be the subset of two lines, which are presented in Fig \ref{q=17-1}. Then $\phi(\mathcal{H})=\mathcal{H}$. It follows that either $\phi(H)=H$, $\phi((\alpha +7)H)=(\alpha +7)H$ or $\phi(H)=(\alpha +7)H$, $\phi((\alpha +7)H)=H$.

If $\phi(H)=H$, $\phi((\alpha +7)H)=(\alpha +7)H$, we have that $aH^{v}=aH=H$ and $a((\alpha+7)H)^{v}=a(\alpha+7)^{v}H=(\alpha+7)H$, then $a\in H$ and $v=1$. So that $\phi(\gamma)=a\gamma$ with $a\in H$.

If $\phi(H)=(\alpha +7)H$, $\phi((\alpha +7)H)=H$, we have that $aH^{v}=aH=(\alpha+7)H$ and $a((\alpha+7)H)^{v}=a(\alpha+7)^{v}H=H$, then $a\in (\alpha+7)H$ and $|v|=2$. So that $\phi(\gamma)=a\gamma^{17}$ with $a\in (\alpha+7)H$.

Set $\sigma(\gamma)=(\alpha+7)\gamma^{17}$. Then $\sigma ^8=1$ and $\lg \sigma\rg=\{\phi_a,\phi_a'\mid \phi_a(\gamma)=a\gamma, \phi_{a'}(\gamma)=a'\gamma^{17}, a\in H, a'\in (\alpha+7)H\}\cong \ZZ_8$. The action of $\sigma$ on the points of clique $C_{17}^{A}$ is presented in the following table. 
 \vskip 1mm
 
 \begin{table}[htp]
\centering
\setlength{\tabcolsep}{0.7mm}
\begin{tabular}{cccccccccc}
\toprule  
$\gamma$ & 0 &$1$ &$4$ &$16$& $13$ & $\alpha+7$& $4(\alpha+7)$& $16(\alpha+7)$& $13(\alpha+7)$\\
	\midrule  
$\sigma(\gamma)$ & 0 &$\alpha+7$ &$4(\alpha+7)$ &$16(\alpha+7)$ & $13(\alpha+7)$& $4$& $16$& $13$& $1$\\

\midrule 
\end{tabular}
\end{table}

Followed by the above arguments, we know that $\sigma (C_{17}^{A})=C_{17}^{A}$. Then $Aut(P(17^2))_{C_{17}^{A}}=\lg \sigma \rg \cong \ZZ_{8}$, where $\sigma(\gamma)=(\alpha+7)\gamma^{17}$ for any $\gamma\in \FF_{17^2}$.

Note that $\left | Aut(P(17^{2})) \right |=\frac{q^2-1}{2}\times q^{2}\times 2=83232$, so we have that $|\mathcal{C}_{17}^{A}|=|Aut(P(17^2)):Aut(P(17^2))_{C_{17}^{A}}|=10404$.
\qed

\vskip 5mm
Set $C_{s,\gamma,i}:=\left \{ s\delta ^{i}x +\gamma   \mid  x\in C_{17}^{A}\cap \FF_{17^2}  \right \} $ where $s\in S_{0}$, $i\in \{ 0,1,2,3\}$, $\gamma \in \FF_{17^{2}}$.

\begin{lem}
Set $\mathcal{C}_{17}^{A}$ be the orbit of $Aut(P(17^2))$ acting on the cliques with $C_{17}^{A}\in \mathcal{C}_{17}$. Then  $\mathcal{C}_{17}^{A}=\{C_{s,\gamma,i}\mid s\in S_{0},i\in \{0,1,2,3\}, \gamma\in\FF_{17^2}\}$.
\end{lem}
\demo
It is obvious that $\{C_{s,\gamma,i}\mid s\in S_{0},i\in \{0,1,2,3\}, \gamma\in\FF_{17^2}\} \subset \mathcal{C}_{17}^{A}$. Now we will prove that $\C_{s,\gamma,i}=C_{s',\gamma',i'}$ if and only if $s=s'$, $\gamma=\gamma'$ and $i=i'$, where $s,s'\in S_{0}$, $i,i'\in\{0,1,2,3\}$ and $\gamma,\gamma'\in\FF_{17^2}$.

Note that the intersection point of two lines in the clique $C_{s,\gamma,i}$ is $\gamma$. If $C_{s,\gamma,i}=C_{s',\gamma',i'}$, then $\gamma=\gamma'$ and the two special lines $\{s\delta^{i}H,s\delta^{i}(\alpha+7)H\}=\{s'\delta^{i'}H,s'\delta^{i'}(\alpha+7)H\}$. 

If $s\delta^{i}H=s'\delta^{i'}H$ and $s\delta^{i}(\alpha+7)H=s'\delta^{i'}(\alpha+7)H$, then $s\delta^i\in s'\delta^{i'}H$. If $s\delta^i= s'\delta^{i'}$, then $s=s'$ and $i=i'$; if $s\delta^i=\delta^4s'\delta^{i'}$, then $\delta^i=\delta^{4+i'}$, a contradiction; if $s\delta^i= \delta^8 s'\delta^{i'}$, then $\delta^i=\delta^{8+i'}$, a contradiction; if $s\delta^i= \delta^{12} s'\delta^{i'}$, then $\delta^i=\delta^{12+i'}$, a contradiction.

If $s\delta^{i}H=s'\delta^{i'}(\alpha+7)H$ and $s'\delta^{i'}H=s\delta^{i}(\alpha+7)H$, then $(\alpha+7)^2\in H$, a contradiction. Now, we have $s=s'$ and $i=i'$.

Followed from the above arguments, we have that  $\mathcal{C}_{17}^{A}=\{C_{s,\gamma,i}\mid s\in S_{0}, i\in \{0,1,2,3\}, \gamma\in\FF_{17^2}\}$, because $|\{C_{s,\gamma,i}\mid s\in S_{0},i\in \{0,1,2,3\}, \gamma\in\FF_{17^2}\}|=\frac{q+1}{2}\times q^{2}\times 4=10404$.
\qed

\subsubsection{$C_{17}^{B}$-construction}
\begin{figure}[htp]
    \centering
    \includegraphics[height=4cm, width=8cm]{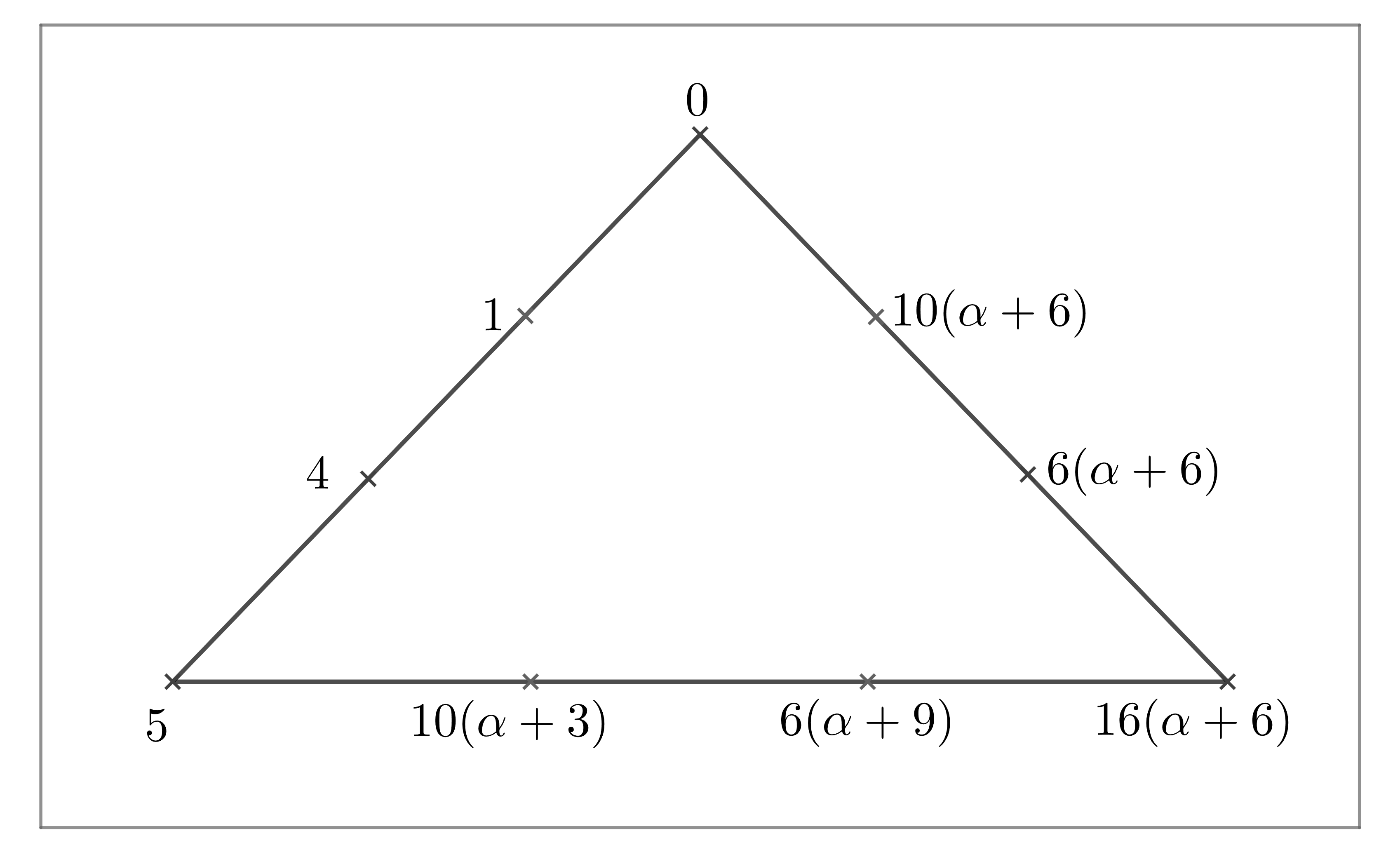}
     \caption{\label{q=17-2}$C_{17}^{B}$}
    \quad
    \end{figure}

Let $H= \{ 1,4,5 \}$  be a subset in $\FF_{17}^{\ast}$. Set $C_{13}^{B}:=\{0\}\cup H\cup10(\alpha+6)H\cup\{10(\alpha+3)\}\cup 
\{6(\alpha+9)\}$ be a subset of finite field $\FF_{17^2}$. By Magma, we have that $C_{17}^{B}$ is a maximal clique in $P(17^2)$ with size $\frac{q+1}{2}$ for $q=17$.
The structure of the clique $C_{17}^{B}$ is presented in Fig \ref{q=17-2}.

 \begin{lem}\label{lem:3.9}
 Set $\mathcal{C}_{17}^{B}$ be the orbit of $Aut(P(17^2))$ acting on the cliques with $C_{17}^{B}\in \mathcal{C}_{17}^{B}$. Then $Aut(P(17^2))_{C_{17}^{B}}=\lg \sigma, \tau  \rg\cong\D_6$, where $ \sigma(\gamma)=10(\alpha+11)\gamma+5$ and $\tau(\gamma)=10(\alpha+6)\gamma^{17}$ for any $\gamma\in \FF_{17^2}$.  Moreover,  $|\mathcal{C}_{17}^{B}|=13872$.    
\end{lem}
\demo
Note that for any $\phi \in \Aut(P(17^2))$,
$\phi(\gamma)=a\gamma^{v}+b$, where $ a\in S$, $b\in \FF_{17^{2}}$, $v\in Gal(\FF_{17^{2}})$.  Set $\mathcal{P}:=\{0,5,16(\alpha+6)\}$ be the three intersecting points in clique $C_{17}^{B}$. Then $\phi(\mathcal{P})=\mathcal{P}$. It follows six cases.

 If $\phi_1(0)=0$, $\phi_1(5)=5$ and $\phi_1(16(\alpha+6))=16(\alpha+6)$, then we have that $a=1$, $b=0$ and $v=1$. It follows that $\phi_1(\gamma)=\gamma$ for any $\gamma\in \FF_{17^2}$.

 If $\phi_2(0)=0$, $\phi_2(5)=16(\alpha+6)$ and $\phi_2(16(\alpha+6))=5$, then we have that $a=10(\alpha+6)$, $b=0$ and $|v|=2$. It follows that $\phi_{2}(\gamma)=10(\alpha+6)\gamma^{17}$ for any $\gamma\in \FF_{17^2}$. 
 
 If $\phi_3(0)=5$, $\phi_3(5)=0$ and $\phi_3(16(\alpha+6))=16(\alpha+6)$, then we have that $a=16$, $b=5$ and $|v|=2$. It follows that $\phi_{3}(\gamma)=16\gamma^{17}+5$ for any $\gamma\in \FF_{17^2}$.

If $\phi_4(0)=5$, $\phi_4(5)=16(\alpha+6)$ and $\phi_4(16(\alpha+6))=0$, then we have that $a=10(\alpha+11)$, $b=5$ and $v=1$. It follows that $\phi_{4}(\gamma)=10(\alpha+11)\gamma+5$ for any $\gamma\in \FF_{17^2}$.

 If $\phi_5(0)=16(\alpha+6)$, $\phi_5(5)=5$ and $\phi_5(16(\alpha+6))=0$, then we have that $a=7(\alpha+11)$, $b=16(\alpha+6)$ and $|v|=2$. It follows that $\phi_{5}(\gamma)=7(\alpha+11)\gamma^{17}+16(\alpha+6)$ for any $\gamma\in \FF_{17^2}$.

 If $\phi_6(0)=16(\alpha+6)$, $\phi_6(5)=0$ and $\phi_6(16(\alpha+6))=5$, then we have that $a=7(\alpha+6)$, $b=16(\alpha+6)$ and $v=1$. It follows that $\phi_{6}(\gamma)=7(\alpha+6)\gamma+16(\alpha+6)$ for any $\gamma\in \FF_{17^2}$.

Set $\sigma=\phi_{4}$, $\tau=\phi_{2}$, then we have that $\sigma^{3}=1$, $\tau^{2}=1$, $\tau^{-1}\sigma\tau=\sigma^{-1}$ and $\lg \sigma, \tau  \rg=\{1,\sigma,\sigma^{2},\tau,\tau\sigma, \tau\sigma^{2}\}=\{\phi_{1},\phi_{4},\phi_{6},\phi_{2}, \phi_{5},\phi_{3}\}$. The action of $\sigma$ and $\tau$ on the points in the clique $C_{17}^{B}$ is presented in the following table. 
 \begin{table}[htp]
\centering
\small
\setlength{\tabcolsep}{0.6mm}
\begin{tabular}{cccccccccc}
\toprule  
$\gamma$ & 0 &$1$ &$4$ &$5$& $10 (\alpha +6   )$ & $6 (\alpha +6 )$& $16 (\alpha +6 )$& $10 (\alpha +3 )$& $6 (\alpha +9 )$\\
	\midrule  
$\sigma(\gamma)$ & $5$ &$10 (\alpha +3 )$ & $6 (\alpha +9 )$ & $16 (\alpha +6 )$& $4$& $1$& $0$& $6 (\alpha +6 )$& $10 (\alpha +6 )$\\
$\tau(\gamma)$ & $0$ &$10 (\alpha +6 )$ & $6 (\alpha +6 )$ & $16 (\alpha +6 )$& $1$& $4$& $5$& $6 (\alpha +9 )$& $10 (\alpha +3 )$\\
\midrule 
\end{tabular}
\end{table}

Followed by the above arguments, we know that $\sigma(C_{17}^{B})=C_{17}^{B}$ and $\tau(C_{17}^{B})=C_{17}^{B}$, then $Aut(P(17^2))_{C_{17}^{B}}=\lg \sigma, \tau  \rg\cong D_{6}$  where $ \sigma(\gamma)=10(\alpha+11)\gamma+5$ and $\tau(\gamma)=10(\alpha+6)\gamma^{17}$ for any $\gamma\in \FF_{17^2}$.

Note that $\left | Aut(P(17^{2})) \right |=\frac{q^2-1}{2}\times q^{2}\times 2=83232$, so we have that $|\mathcal{C}_{17}^{B}|=|Aut(P(17^2)):Aut(P(17^2))_{C_{17}^{B}}|=13872$.
\qed

\vskip 5mm
Set $C_{s,\gamma}:=\left \{ sx +\gamma   \mid  x\in C_{17}^{B}\cap \FF_{17^2}  \right \} $ where $s\in S$, $\gamma \in \FF_{17^{2}}$.

\begin{lem}
Set $\mathcal{C}_{17}^{B}$ be the orbit of $Aut(P(17^2))$ acting on the cliques with $C_{17}^{B}\in \mathcal{C}_{17}$. Then  $\mathcal{C}_{17}^{B}=\{C_{s,\gamma}\mid s\in S,\gamma\in\FF_{17^2}\}$.
\end{lem}

\demo
It is obvious that $\{C_{s,\gamma}\mid s\in S, \gamma\in\FF_{17^2}\} \subset \mathcal{C}_{17}^{B}$. If $\gamma=\gamma'$, now we will prove that $\C_{s,\gamma}=C_{s',\gamma'}$ if and only if $s=s'$,  where $s,s'\in S$, $\gamma,\gamma'\in\FF_{17^2}$.

 If $C_{s,\gamma}=C_{s',\gamma'}$, then the three special points $\{\gamma,5s+\gamma,16(\alpha+6)s+\gamma\}=\{\gamma',5s'+\gamma',16(\alpha+6)s'+\gamma'\}$. If $5s+\gamma=16(\alpha+6)s'+\gamma'$ and $5s'+\gamma'=16(\alpha+6)s+\gamma$, then $(\alpha+6)=12$, a contradiction. Now, we have $s=s'$.
From Lemma \ref{lem:3.9}, we have that $C_{17}^{B}=\phi _{4}(C_{17}^{B})=\phi_{6}(C_{17}^{B})$, where $\lg \phi_4,\phi_6\rg\cong \ZZ_3$ permutate the three intersecting points. 
And then $\mathcal{C}_{17}^{B}=\{C_{s,\gamma}\mid s\in S, \gamma\in\FF_{17^2}\}$, because $|\{C_{s,\gamma}\mid s\in S,\gamma\in\FF_{17^2}\}|=\frac{q^{2}-1}{2}\times q^{2}\times \frac{1}{3}=13872$.
\qed

\subsubsection{$C_{17}^{C}$-construction}
Let $H=\left \{ \delta^{4} ,\delta ^{8},\delta ^{12} , \delta ^{16} \right \}=\left \{ 4,16,13,1 \right \} $ be a subgroup in $\FF_{17}^{\ast}$ with $\delta:=10$. Then
$\FF_{17}^{\ast}=H \cup \delta H \cup \delta ^{2} H \cup \delta ^{3} H=\left \{ 4,16,13,1 \right \}\cup  \left \{6,7,11,10  \right \}\cup  \left \{ 9,2,8,5\right \}\cup\left \{ 5,3,12,14  \right \}$.

\begin{figure}[htp]
    \centering
    \includegraphics[height=4cm, width=8cm]{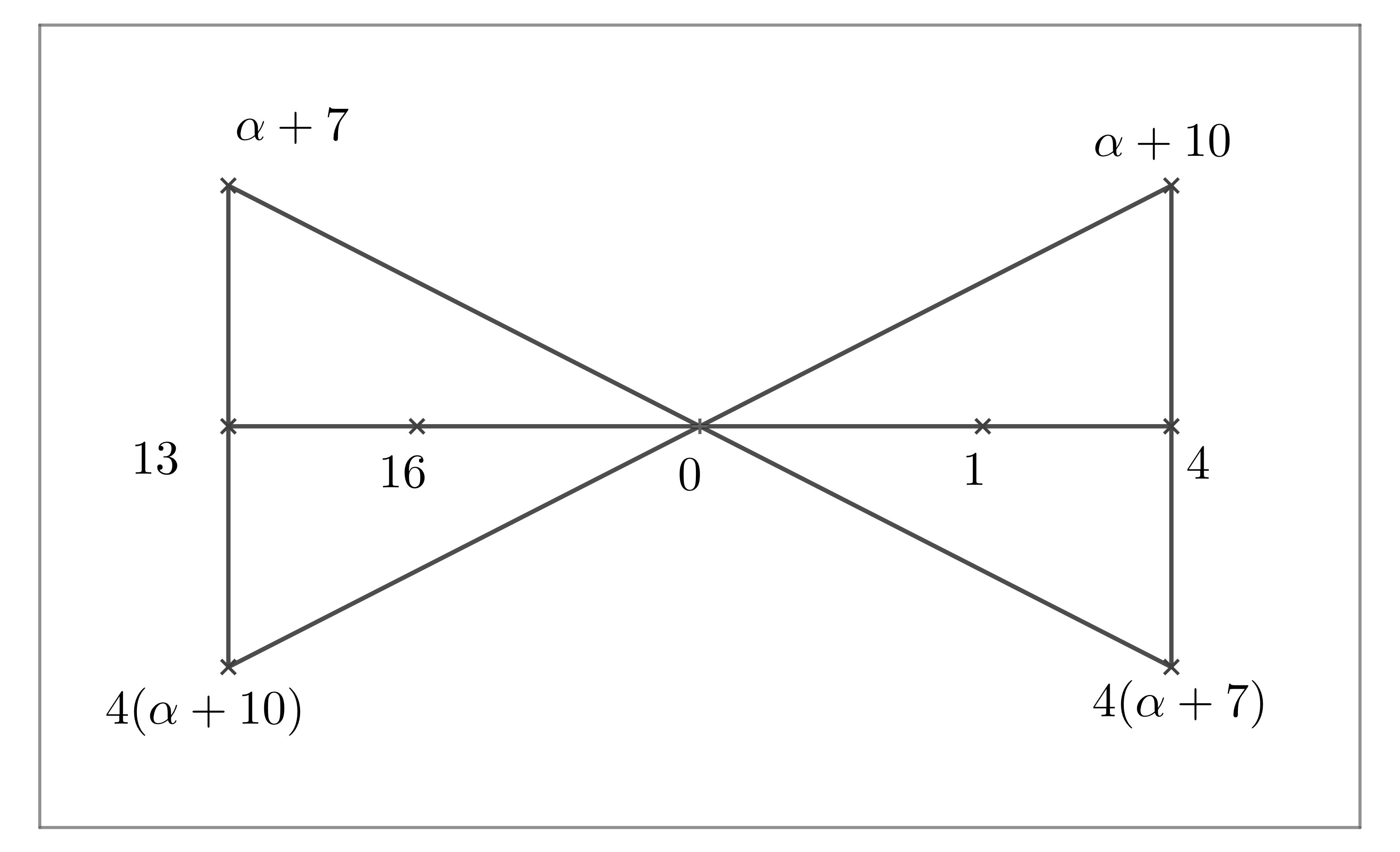}
    \quad
     \caption{\label{q=17-3}$C_{17}^{C}$}
    \end{figure}

Set $C_{17}^{C}:=\{0\}\cup H\cup (\alpha+7)\{1,4\}\cup (\alpha+10)\{1,4\}$ be a subset of finite field $\FF_{17^2}$. By Magma, we have that $C_{17}^{C}$ is a maximal clique in $P(17^2)$ with size $\frac{q+1}{2}$ for $q=17$.
The structure of the clique $C_{17}^{C}$ is presented in Fig \ref{q=17-3}.

\begin{lem}
Set $\mathcal{C}_{17}^{C}$ be the orbit of $Aut(P(17^2))$ acting on the cliques with $C_{17}^{C}\in \mathcal{C}_{17}^{C}$. Then $Aut(P(17^2))_{C_{17}^{C}}=\lg \phi' \rg\cong \ZZ_{2}$, where $\phi'(\gamma)=16\gamma^{17}$ for any $\gamma\in \FF_{17^2}$. Moreover, $|\mathcal{C}_{17}^{C}|=41616$.
\end{lem}
\demo
Note that for any $\phi \in \Aut(P(17^2))$,
$\phi(\gamma)=a\gamma^{v}+b$, where $ a\in S$, $b\in \FF_{17^{2}}$, $v\in Gal(\FF_{17^{2}})$. 
Note that the point $0$ in the clique $C_{17}^C$ is a intersecting point of three lines, and one of these lines contains five points in this clique. 
If $\phi\in Aut(P(17^{2}))_{C_{17}^{C}}$, then  $\phi(0)=0$ and $b=0$. And $\phi(H)=H$, we have that $aH^{v}=aH=H$, then $a\in H$. So that $\phi(\gamma)=a\gamma^{v}$ with $a\in H$, $v\in Gal(\FF_{17^{2}})$. Set $H'=\{1,4\}$ be a subset in $\FF_{q}^{\ast}$ and $\mathcal{L}:=\{H'(\alpha+7), H'(\alpha+10)\}$ be the subset of two lines, which are presented in Fig \ref{q=17-3}. Then $\phi(\mathcal{L})=\mathcal{L}$. It follows that either $\phi(H'(\alpha+7))=H'(\alpha+7)$, $\phi(H'(\alpha+10))=H'(\alpha+10)$ or $\phi(H'(\alpha+7))=H'(\alpha+10)$, $\phi(H'(\alpha+10))=H'(\alpha+7)$.

If $\phi(H'(\alpha+7))=H'(\alpha+7)$ 
 and $\phi(H'(\alpha+10))=H'(\alpha+10)$, then we have that $a(H'(\alpha+7))^{v}=aH'(\alpha+7)^{v}=H'(\alpha+7)$ and $a(H'(\alpha+10))^{v}=aH'(\alpha+10)^{v}=H'(\alpha+10)$, then $a=1$ and $v=1$. So  $\phi(\gamma)=\gamma$  for any $\gamma\in \FF_{17^2}$.

If $\phi'(H'(\alpha+7))=H'(\alpha+10)$
 and $\phi'(H'(\alpha+10))=H'(\alpha+7)$, then we have that $a(H'(\alpha+7))^{v}=aH'(\alpha+7)^{v}=H'(\alpha+10)$ and $a(H'(\alpha+10))^{v}=aH'(\alpha+10)^{v}=H'(\alpha+7)$, then $a=16$ and $|v|=2$. So  $\phi'(\gamma)=16\gamma^{17}$ for any $\gamma\in \FF_{17^2}$. 

 The action of $\phi'$ acting on the points in the clique $C_{17}^{C}$ is presented by the following table. 
 \begin{table}[htp]
\centering
\begin{tabular}{cccccccccc}
\toprule  
$\gamma$ & 0 &$1$ &$4$ &$16$& $13$ & $\alpha+7$& $4(\alpha+7)$& $\alpha+10$& $4(\alpha+10)$\\
	\midrule  
$\phi'(\gamma)$ & 0 &$16$ &$13$ &$1$ & $4$& $\alpha+10$& $4(\alpha+10)$& $\alpha+7$& $4(\alpha+7)$\\
\midrule 
\end{tabular}
\end{table}

Followed by the above arguments, we know that $\phi'(C_{17}^{C})=C_{17}^{C}$ and $\phi'^2=1$. Then $Aut(P(17^2))_{C_{17}^{C}}=\lg \phi'  \rg\cong \ZZ_{2}$  where $ \phi'(\gamma)=16\gamma^{17}$ for any $\gamma\in \FF_{17}^2$.

Note that $\left | Aut(P(17^{2})) \right |=\frac{q^2-1}{2}\times q^{2}\times 2=83232$, so we have that $|\mathcal{C}_{17}^{C}|=|Aut(P(17^2)):Aut(P(17^2))_{C_{17}^{C}}|=41616$.
\qed

\vskip 5mm
Set $C_{s,\gamma}:=\left \{ sx +\gamma   \mid  x\in C_{17}^{C}\cap \FF_{17^2}  \right \} $ where $s\in S$, $\gamma \in \FF_{17^{2}}$.

\begin{lem}
Set $\mathcal{C}_{17}^{C}$ be the orbit of $Aut(P(17^2))$ acting on the cliques with $C_{17}^{C}\in \mathcal{C}_{17}$. Then  $\mathcal{C}_{17}^{C}=\{C_{s,\gamma}\mid s\in S,\gamma\in\FF_{17^2}\}$.
\end{lem}
\demo
It is obvious that $\{C_{s,\gamma}\mid s\in S, \gamma\in\FF_{17^2}\} \subset \mathcal{C}_{17}^{C}$. Now we will prove that $\C_{s,\gamma}=C_{s',\gamma'}$ if and only if $s=s'$ and $\gamma=\gamma'$, where $s,s'\in S$ and $\gamma,\gamma'\in\FF_{17^2}$.

Note that the intersecting point of three special lines in the clique $C_{s,\gamma}$ is $\gamma$. If $C_{s,\gamma}=C_{s',\gamma'}$, then $\gamma=\gamma'$ and the subset of two lines $\{s(\alpha+7)H'+\gamma,s(\alpha+10)H'+\gamma\}=\{s'(\alpha+7)H'+\gamma',s'(\alpha+10)H'+\gamma'\}$. If $s(\alpha+7)H'+\gamma=s'(\alpha+7)H'+\gamma'$ and $s(\alpha+10)H'+\gamma=s'(\alpha+10)H'+\gamma'$, then $s=s'$;   if $s(\alpha+7)H'+\gamma=s'(\alpha+10)H'+\gamma'$ and $s(\alpha+10)H'+\gamma=s'(\alpha+7)H'+\gamma'$, then $(\alpha+7)^2H'=(\alpha+10)^2H'$, a contradiction. Now, we have $s=s'$. 

And then $\mathcal{C}_{17}^{C}=\{C_{s,\gamma}\mid s\in S, \gamma\in\FF_{17^2}\}$, because $|\{C_{s,\gamma}\mid s\in S, \gamma\in\FF_{17^2}\}|=\frac{q^{2}-1}{2}\times q^{2}=41616$.
\qed

\subsubsection{$C_{17}^{D}$-construction}

Let $H=\left \{ \delta^{4} ,\delta ^{8},\delta ^{12} , \delta ^{16} \right \}=\left \{ 4,16,13,1 \right \} $ be a subgroup in $\FF_{17}^{\ast}$ with $\delta:=10$. Then
$\FF_{17}^{\ast}=H \cup \delta H \cup \delta ^{2} H \cup \delta ^{3} H=\left \{ 4,16,13,1 \right \}\cup  \left \{6,7,11,10  \right \}\cup  \left \{ 9,2,8,5\right \}\cup\left \{ 5,3,12,14  \right \}$.

\begin{figure}[htp]
    \centering
    \includegraphics[height=4cm, width=8cm]{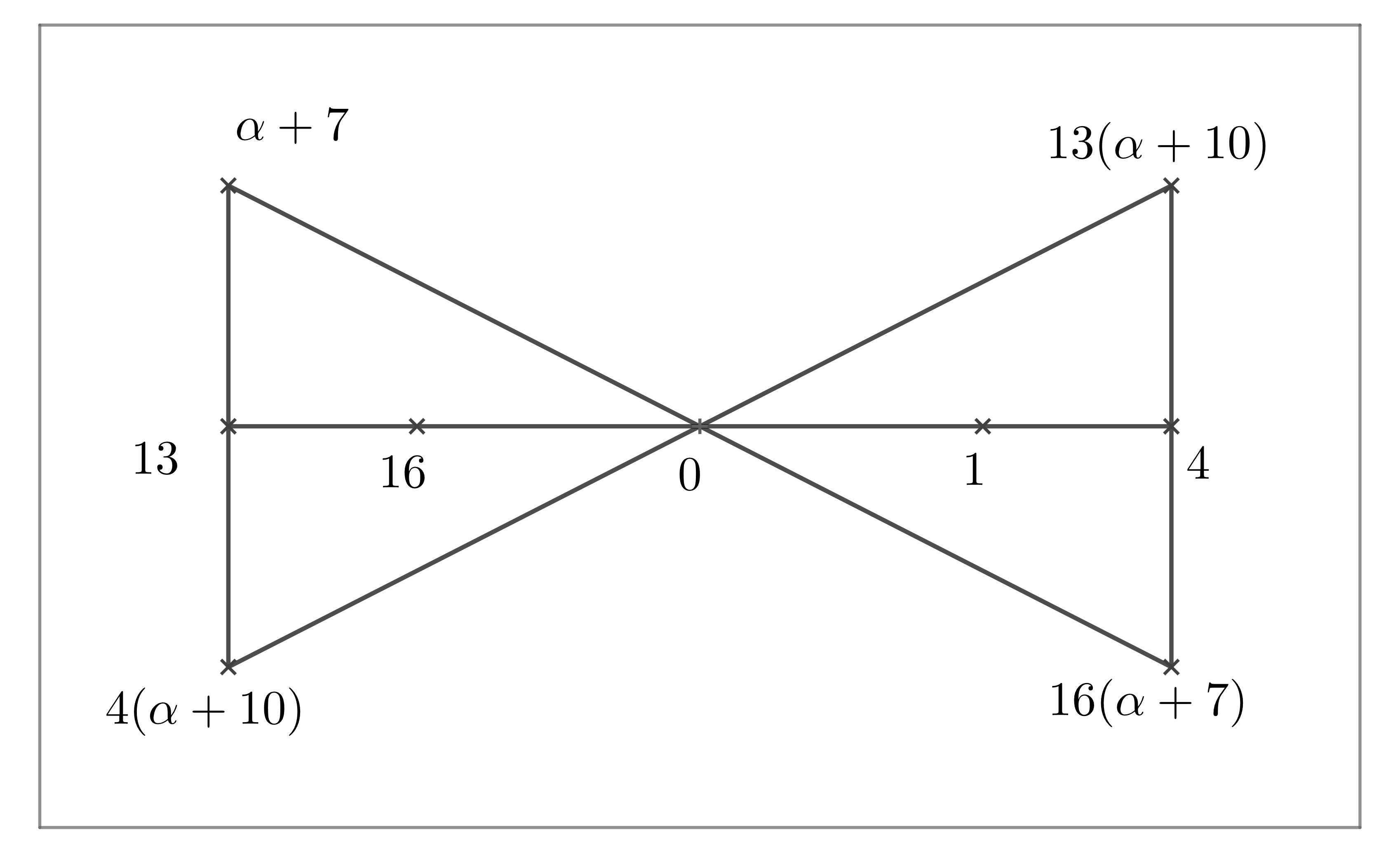}
    \quad
     \caption{\label{q=17-4}$C_{17}^{D}$}
Remarks: There exists other lines that contain at least three points in $C_{17}^{D}$. They are: $\{1, \alpha+7,13(\alpha+10)\}$,$\{16, 4(\alpha+10),16(\alpha+7)\}$.
    \end{figure}

Set $C_{17}^{D}:=\{0\}\cup H\cup (\alpha+7)\{1,16\}\cup (\alpha+10)\{4,13\}$ be a subset of finite field $\FF_{17^2}$. By Magma, we have that $C_{17}^{D}$ is a maximal clique in $P(17^2)$ with size $\frac{q+1}{2}$ for $q=17$.
The structure of the clique $C_{17}^{D}$ is presented in Fig \ref{q=17-4}.

\begin{lem}
 Set $\mathcal{C}_{17}^{D}$ be the orbit of $Aut(P(17^2))$ acting on the cliques with $C_{17}^{D}\in \mathcal{C}_{17}^{D}$. Then $Aut(P(17^2))_{C_{17}^{D}}=\lg \sigma\rg\cong \ZZ_{4}$, where $\sigma(\gamma)=4\gamma^{17}$ for any $\gamma\in \FF_{17^2}$. Moreover, $|\mathcal{C}_{17}^{D}|=20808$.
\end{lem}

\demo
Note that for any $\phi \in \Aut(P(17^2))$,
$\phi(\gamma)=a\gamma^{v}+b$, where $ a\in S$, $b\in \FF_{17^{2}}$, $v\in Gal(\FF_{{17}^2})$.
Note that the point $0$ in the clique $C_{17}^D$ is a intersecting point of three lines, and one of these lines contains five points in this clique. 
If $\phi\in Aut(P(17^{2}))_{C_{17}^{D}}$, then  $\phi(0)=0$ and $b=0$. And $\phi(H)=H$, we have that $aH^{v}=aH=H$, then $a\in H$. So that $\phi(\gamma)=a\gamma^{v}$ with $a\in H$, $v\in Gal(\FF_{{17}^{2}})$. Let $H'=\{1,16\}$ be a subset in $\FF_{17}^{\ast}$ and $H'=-H'$. Set $\mathcal{L}:=\{H'(\alpha+7), 4H'(\alpha+10)\}$ be the subset of two lines, which are presented in Fig \ref{q=17-4}. Then $\phi(\mathcal{L})=\mathcal{L}$. It follows that either $\phi(H'(\alpha+7))=H'(\alpha+7)$, $\phi(4H'(\alpha+10))=4H'(\alpha+10)$ or $\phi(H'(\alpha+7))=4H'(\alpha+10)$, $\phi(4H'(\alpha+10))=H'(\alpha+7)$.

If $\phi(H'(\alpha+7))=H'(\alpha+7)$
 and $\phi(4H'(\alpha+10))=4H'(\alpha+10)$ , then we have that $a(H'(\alpha+7))^{v}=aH'(\alpha+7)^{v}=H'(\alpha+7)$ and $a(4H'(\alpha+10))^{v}=4aH'(\alpha+10)^{v}=4H'(\alpha+10)$, then $a\in\{1,-1\}$ and $v=1$. Now set $\phi_{1}(\gamma)=\gamma$ and $\phi_{2}(\gamma)=-\gamma$ for any $\gamma\in \FF_{17^2}$. 

If $\phi(H'(\alpha+7))=4H'(\alpha+10)$
 and $\phi(4H'(\alpha+10))=H'(\alpha+7)$, then we have that $a(H'(\alpha+7))^{v}=aH'(\alpha+7)^{v}=4H'(\alpha+10)$ and $a(4H'(\alpha+10))^{v}=4aH'(\alpha+10)^{v}=H'(\alpha+7)$, then $a\in\{4,-4\}$ and $|v|=2$. Now set $\phi_{3}(\gamma)=4\gamma^{17}$, $\phi_{4}(\gamma)=-4\gamma^{17}$ for any $\gamma\in \FF_{17^2}$. 

 Set $\sigma=\phi_{3}$, now we have that $\sigma^4=1$ and $\lg \sigma \rg=\{1,\sigma,\sigma^2,\sigma^3 \}=\{\phi_{1},\phi_{3},\phi_{2},\phi_{4}\}\cong\ZZ_4$. The action of $\sigma$ on the points in clique $C_{17}^{D}$ is presented in the following table. 
 
 \begin{table}[htp]
\centering
\begin{tabular}{cccccccccc}
\toprule  
$\gamma$ & 0 &$1$ &$4$ &$16$& $13$ & $\alpha+7$& $16(\alpha+7)$& $4(\alpha+10)$& $13(\alpha+10)$\\
	\midrule  
$\sigma(\gamma)$ & 0 &$4$ &$16$ &$13$ & $1$& $13(\alpha+10)$& $4(\alpha+10)$& $\alpha+7$& $16(\alpha+7)$\\
\midrule 
\end{tabular}
\end{table}

Followed by the above arguments, we know that $\sigma(C_{17}^{D})=C_{17}^{D}$ and then $Aut(P(17^2))_{C_{17}^{D}}=\lg \sigma \rg\cong \ZZ_{4}$  where $ \sigma(\gamma)=4\gamma^{17}$ for any $\gamma\in \FF_{{17}^2}$.

Note that $\left | Aut(P(17^{2})) \right |=\frac{q^2-1}{2}\times q^{2}\times 2=83232$, so we have that $|\mathcal{C}_{17}^{D}|=|Aut(P(17^2)):Aut(P(17^2))_{C_{17}^{D}}|=20808$.
\qed

\vskip 5mm
Set $C_{s,\gamma,i}:= \{ s\delta ^{i}x +\gamma   \mid  x\in C_{17}^{D}\cap \FF_{17^2}  \} $ where $s\in S_{0}$, $i\in \{ 1\dots8\}$, $\gamma \in \FF_{17^{2}}$. 

\begin{lem}
Set $\mathcal{C}_{17}^{D}$ be the orbit of $Aut(P(17^2))$ acting on the cliques with $C_{17}^{D}\in \mathcal{C}_{17}$. Then  $\mathcal{C}_{17}^{D}=\{C_{s,\gamma,i}\mid s\in S_{0},i\in \{1\dots 8\}, \gamma\in\FF_{17^2}\}$.
\end{lem}

\demo
It is obvious that $\{C_{s,\gamma,i}\mid s\in S_{0},i\in \{1\dots8\}, \gamma\in\FF_{17^2}\} \subset \mathcal{C}_{17}^{D}$. Now we will prove that $\C_{s,\gamma,i}=C_{s',\gamma',i'}$ if and only if $s=s'$, $\gamma=\gamma'$ and $i=i'$, where $s,s'\in S_{0}$, $i,i'\in\{1\dots 8\}$ and $\gamma,\gamma'\in\FF_{17^2}$.

Note that the intersection point of three lines in the clique $C_{s,\gamma,i}$ is $\gamma$. If $C_{s,\gamma,i}=C_{s',\gamma',i'}$, then $\gamma=\gamma'$ and $\{s\delta^{i}, 4s\delta^{i}, 13s\delta^{i}, 16s\delta^{i}\}=\{s'\delta^{i'}, 4s'\delta^{i'}, 13s'\delta^{i'}, 16s'\delta^{i'}\}$. Now we have that $s\delta^{i}\in \{s'\delta^{i'}, 4s'\delta^{i'}, 13s'\delta^{i'}, 16s'\delta^{i'}\}$. It follows that $s=s'$ and $i=i'$.

And then $\mathcal{C}_{17}^{D}=\{C_{s,\gamma,i}\mid s\in S_{0}, i\in \{1\dots8\}, \gamma\in\FF_{17^2}\}$, because $|\{C_{s,\gamma,i}\mid s\in S_{0},i\in \{1\dots8\}, \gamma\in\FF_q^2\}|=\frac{q+1}{2}\times q^{2}\times 8=20808$.
\qed

\subsubsection{$C_{17}^{E}$-construction }
Let $H=\left \{ \delta^{4} ,\delta ^{8},\delta ^{12} , \delta ^{16} \right \}=\left \{ 4,16,13,1 \right \} $ be a subgroup in $\FF_{17}^{\ast}$ with $\delta:=10$. Then
$\FF_{17}^{\ast}=H \cup \delta H \cup \delta ^{2} H \cup \delta ^{3} H=\left \{ 4,16,13,1 \right \}\cup  \left \{6,7,11,10  \right \}\cup  \left \{ 9,2,8,5\right \}\cup\left \{ 5,3,12,14  \right \}$.

\begin{figure}[htp]
    \centering
    \includegraphics[height=4cm, width=8cm]{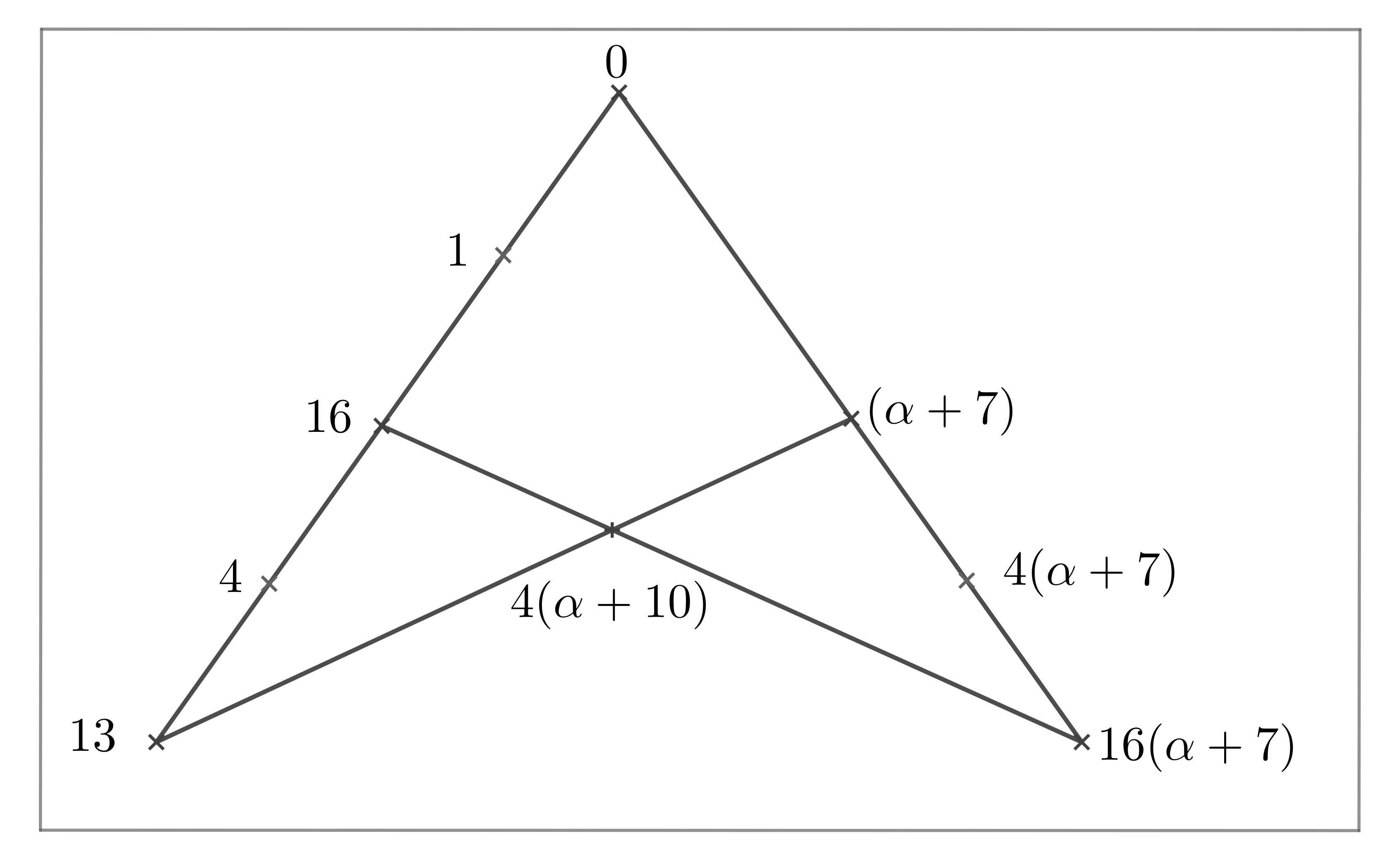}
    \quad
     \caption{\label{q=17-5}$C_{17}^{E}$}
    \end{figure}

Set $C_{17}^{E}:=\{0\}\cup H\cup (\alpha+7)\{1,4,16\}\cup\{4(\alpha+10)\}$ be a subset of finite field  $\FF_{17^2}$. By Magma, we have that $C_{17}^{E}$ is a maximal clique in $P(17^2)$ with size $\frac{q+1}{2}$ for $q=17$.
The structure of the clique $C_{17}^{E}$ is presented in Fig \ref{q=17-5}.

\begin{lem}
Set $\mathcal{C}_{17}^{E}$ be the orbit of $Aut(P(17^2))$ acting on the cliques with $C_{17}^{E}\in \mathcal{C}_{17}^{E}$. Then $Aut(P(17^2))_{C_{17}^{E}}=1$ and $|\mathcal{C}_{17}^{E}|=83232$. 
\end{lem}
\demo
Note that for any $\phi \in \Aut(P(17^2))$,
$\phi(\gamma)=a\gamma^{v}+b$, where $ a\in S$, $b\in \FF_{17^{2}}$, $v\in Gal(\FF_{17^{2}})$.
Note that the point $0$ in the clique $C_{17}^E$ is a intersecting point of two lines, where one line contain five points in this clique and the other line contain four points; the point $4(\alpha+10)$ in the clique $C_{17}^E$ is a intersecting point of two lines, which contain three points separately in this clique. 
If $\phi\in Aut(P(17^{2}))_{C_{17}^{E}}$, then  $\phi(0)=0$, $\phi(4(\alpha+10))=4(\alpha+10)$. So we have that $b=0$ and $a(4(\alpha+10))^v=4a(\alpha+10)^v=4(\alpha+10)$, then $a=1 $ and $v=1$. Followed by the above arguments, we know that $|Aut(P(17^2))_{C_{17}^{E}}|= 1 $.

Note that $\left | Aut(P(17^{2})) \right |=\frac{q^2-1}{2}\times q^{2}\times 2=83232$, so we have that $|\mathcal{C}_{17}^{E}|=|Aut(P(17^2)):Aut(P(17^2))_{C_{17}^{E}}|=83232$.
\qed
\vskip 5mm

Set $C_{s,\gamma}^{\eta}:=\{sx^{\eta}+\gamma| x\in C_{17}^{E}\cap \FF_{17^2}\}$ for $\eta\in\{1,17\}$, where $s\in S$ and $\gamma\in\FF_{17^2}$ . 

\begin{lem}
Set $\mathcal{C}_{17}^{E}$ be the orbit of $Aut(P(17^2))$ acting on the cliques with $C_{17}^{E}\in \mathcal{C}_{17}^{E}$. Then $\mathcal{C}_{17}^{E}=\{C_{s,\gamma}\mid s\in S, \gamma\in\FF_{17^2}\}\cup \{C_{s,\gamma}^{17}\mid s\in S, \gamma\in\FF_{17^2}\}$.
\end{lem}
\demo
It is obvious that $\{\{C_{s,\gamma}\mid s\in S, \gamma\in\FF_{17^{2}}\}\cup \{C_{s,\gamma}^{17}\mid s\in S, \gamma\in\FF_{17^2}\}\}\subset \mathcal{C}_{17}^{E}$. Now we will prove that $C_{s,\gamma}^{\eta}=C_{s',\gamma'}^{\eta'}$ if and only if $s=s'$, $\gamma=\gamma'$ and $\eta=\eta'$, where $s,s'\in S$, $\gamma,\gamma'\in\FF_{17^{2}}$ and $\eta,\eta'\in\{1,17\}$.
 If $C_{s,\gamma}^{\eta}=C_{s',\gamma'}^{\eta'}$, then $\gamma=\gamma'$, the subset $sH^{\eta} =s'H^{\eta'}$ and  $4(\alpha+10)^{\eta}s=4(\alpha+10)^{\eta'}s'$. Then $s\in\{s',4s',16s',13s'\}$, and we have $s=s'$ and $\eta=\eta'$.

And then $\mathcal{C}_{17}^{E}=\{C_{s,\gamma}\mid s\in S, \gamma\in\FF_{17^2}\}\cup \{C_{s,\gamma}^{17}\mid s\in S, \gamma\in\FF_{17^2}\}$, because $|\{C_{s,\gamma}\mid s\in S, \gamma\in\FF_{17^2}\}\cup \{C_{s,\gamma}^{17}\mid s\in S, \gamma\in\FF_{17^2}\}|=\frac{q^2-1}{2}\times q^{2}\times 2=83232$.
\qed

\subsubsection{$C_{17}^{F}$-construction}

\begin{figure}[htp]
    \centering
    \includegraphics[height=4cm, width=8cm]{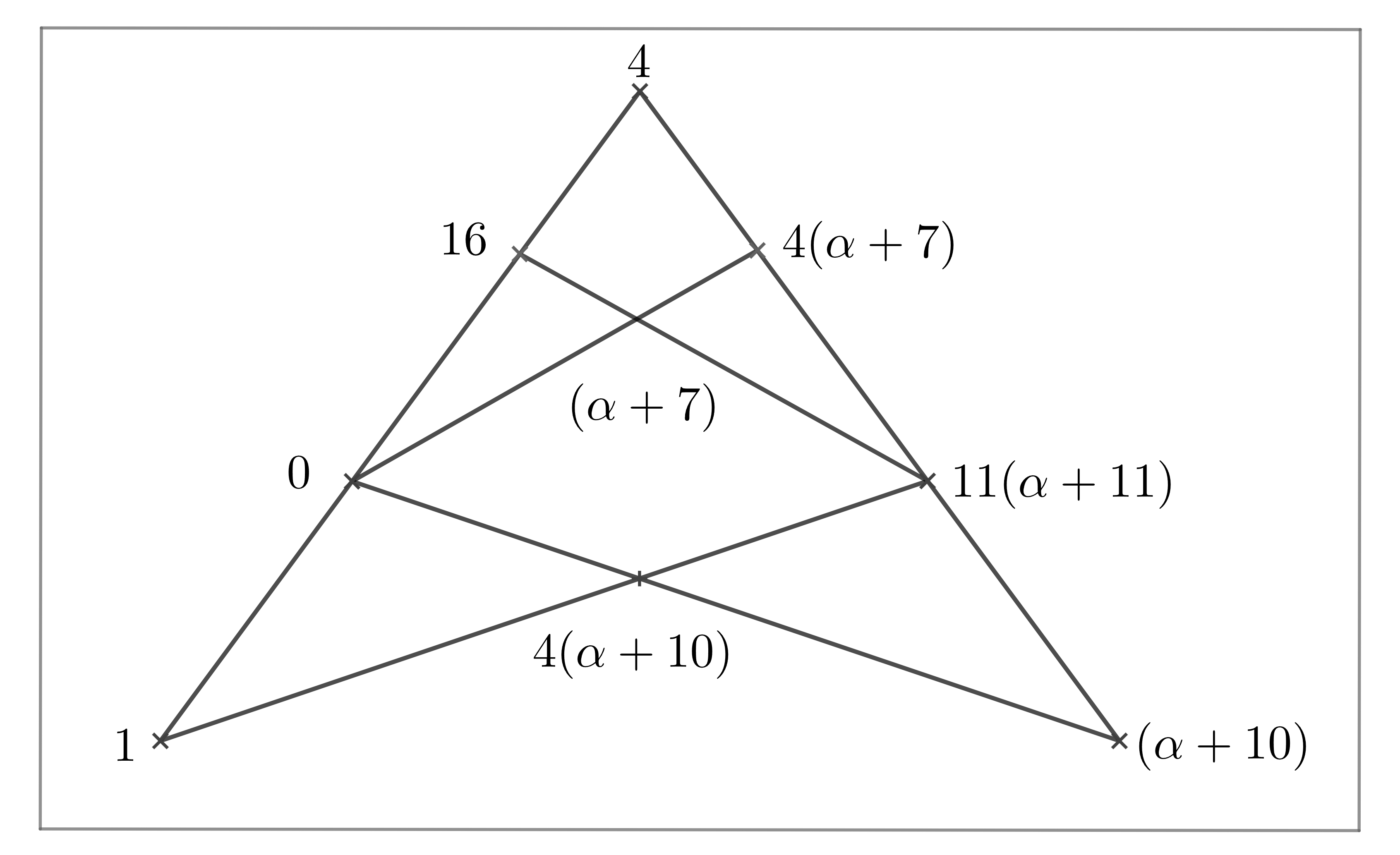}
    \quad
     \caption{\label{q=17-6}$C_{17}^{F}$}
    \end{figure}

 Let $H= \{ 1,4,16\}$  be a subset in $\FF_{q}^{\ast}$. Set $C_{17}^{F}:=\{0\}\cup H\cup (\alpha+7)\{1,4\}\cup (\alpha+10)\{1,4\}\cup \{11(\alpha+11)\}$ be a subset of finite field  $\FF_{17^2}$. By Magma, we have that $C_{17}^{F}$ is a maximal clique in $P(17^2)$ with size $\frac{q+1}{2}$ for $q=17$.
The structure of the clique $C_{17}^{F}$ is presented in Fig \ref{q=17-6}.

\begin{lem}
Set $\mathcal{C}_{17}^{F}$ be the orbit of $Aut(P(17^2))$ acting on the cliques with $C_{17}^{F}\in \mathcal{C}_{17}^{F}$. Then $Aut(P(17^2))_{C_{17}^{F}}=\lg \phi'  \rg\cong \ZZ_{2}$  where $\phi'(\gamma)=10(\alpha+6)\gamma^{17}+11(\alpha+11)$ for any $\gamma\in \FF_{17^2}$. Moreover, $|\mathcal{C}_{17}^{F}|=41616$. 
\end{lem}
\demo
Note that for any $\phi \in \Aut(P(17^2))$,
$\phi(\gamma)=a\gamma^{v}+b$, where $ a\in S$, $b\in \FF_{17^{2}}$, $v\in Gal(\FF_{17^{2}})$. If $\phi\in Aut(P(17^{2}))_{C_{17}^{F}}$, then $\phi(4)=4$ and $4a+b=4$. Set $\mathcal{P}:=\{0,11(\alpha+11)\}$ be the two special  points which are the intersecting points of three lines in  clique $C_{17}^{F}$. Then $\phi(\mathcal{P})=\mathcal{P}$. It follows that either $\phi(0)=0$, $\phi(11(\alpha+11))=11(\alpha+11)$ or $\phi(0)=11(\alpha+11)$, $\phi(11(\alpha+11))=0$.

If  $\phi(0)=0$ and $\phi(11(\alpha+11))=11(\alpha+11)$, then we have that $a=1$, $b=0$ and $v=1$. So that $\phi(\gamma)=\gamma$ for any $\gamma\in \FF_{17^2}$.

If $\phi'(0)=11(\alpha+11)$ and $\phi'(11(\alpha+11))=0$, then we have that $a=10(\alpha+6)$, $b=11(\alpha+11)$ and $|v|=2$. So that $\phi'(\gamma)=10(\alpha+6)\gamma^{17}+11(\alpha+11)$ and $\phi'^2=1$, $\lg\phi'\rg=\{\phi,\phi'\}\cong \ZZ_2$

The action of $\phi'$ acting on the points in the clique $C_{17}^{F}$ is presented by the following table. 
 \begin{table}[htp]
\centering
\setlength{\tabcolsep}{0.6mm}
\begin{tabular}{cccccccccc}
\toprule  
$\gamma$ & 0 &$1$ &$4$ &$16$ & $\alpha+7$& $4(\alpha+7)$& $\alpha+10$& $4(\alpha+10)$& $11(\alpha+11)$\\
	\midrule  
$\phi'(\gamma)$ & $11(\alpha+11)$ &$4(\alpha+7)$ &$4$ &$\alpha+10$ & $4(\alpha+10)$& $1$& $16$& $\alpha+7$& $0$\\
\midrule 
\end{tabular}
\end{table}

Followed by the above arguments, we know that $\phi'(C_{17}^{F})=C_{17}^{F}$ then $Aut(P(17^2))_{C_{17}^{F}}=\lg \phi'  \rg\cong \ZZ_{2}$  where $ \phi'(\gamma)=10(\alpha+6)\gamma^{17}+11(\alpha+11)$ with $\gamma\in \FF_{{17}^2}$.

Note that $\left | Aut(P(17^{2})) \right |=\frac{q^2-1}{2}\times q^{2}\times 2=83232$, so we have that $|\mathcal{C}_{17}^{F}|=|Aut(P(17^2)):Aut(P(17^2))_{C_{17}^{F}}|=41616$.
\qed
\vskip 5mm
Set $C_{s,\gamma}:=\left \{ sx +\gamma   \mid  x\in C_{17}^{F}\cap \FF_{17^2}  \right \} $ where $s\in S$, $\gamma \in \FF_{q^{2}}$. 

\begin{lem}
Set $\mathcal{C}_{17}^{F}$ be the orbit of $Aut(P(17^2))$ acting on the cliques with $C_{17}^{F}\in \mathcal{C}_{17}$. Then  $\mathcal{C}_{17}^{F}=\{C_{s,\gamma}\mid s\in S, \gamma\in\FF_{17^2}\}$.
\end{lem}

\demo
It is obvious that $\{C_{s,\gamma}\mid s\in S, \gamma\in\FF_{17^2}\} \subset \mathcal{C}_{17}^{F}$. Now we will prove that $\C_{s,\gamma}=C_{s',\gamma'}$ if and only if $s=s'$, $\gamma=\gamma'$  where $s,s'\in S$ and $\gamma,\gamma'\in\FF_{17^2}$.

If $C_{s,\gamma}=C_{s',\gamma'}$, then $4s+\gamma=4s'+\gamma'$ and the two special points $\{\gamma, 11(\alpha+11)s+\gamma \}=\{\gamma', 11(\alpha+11)s'+\gamma'\}$. If $\gamma=\gamma'$ and $11(\alpha+11)s+\gamma=11(\alpha+11)s'+\gamma'$, then $s=s'$ and $\gamma=\gamma'$. If $\gamma=11(\alpha+11)s'+\gamma'$ and $11(\alpha+11)s+\gamma=\gamma'$, then $(\alpha+11)=10$, a contradiction. Now, we have $s=s'$ and $\gamma=\gamma'$.

And then $\mathcal{C}_{17}^{F}=\{C_{s,\gamma}\mid s\in S,\gamma\in\FF_{17^2}\}$, because $|\{C_{s,\gamma}\mid s\in S, \gamma\in\FF_{17^2}\}|=\frac{q^{2}-1}{2}\times q^{2}=41616$.
\qed
\subsubsection{$C_{17}^{G}$-construction}
 
\begin{figure}[htp]
    \centering
    \includegraphics[height=4cm, width=8cm]{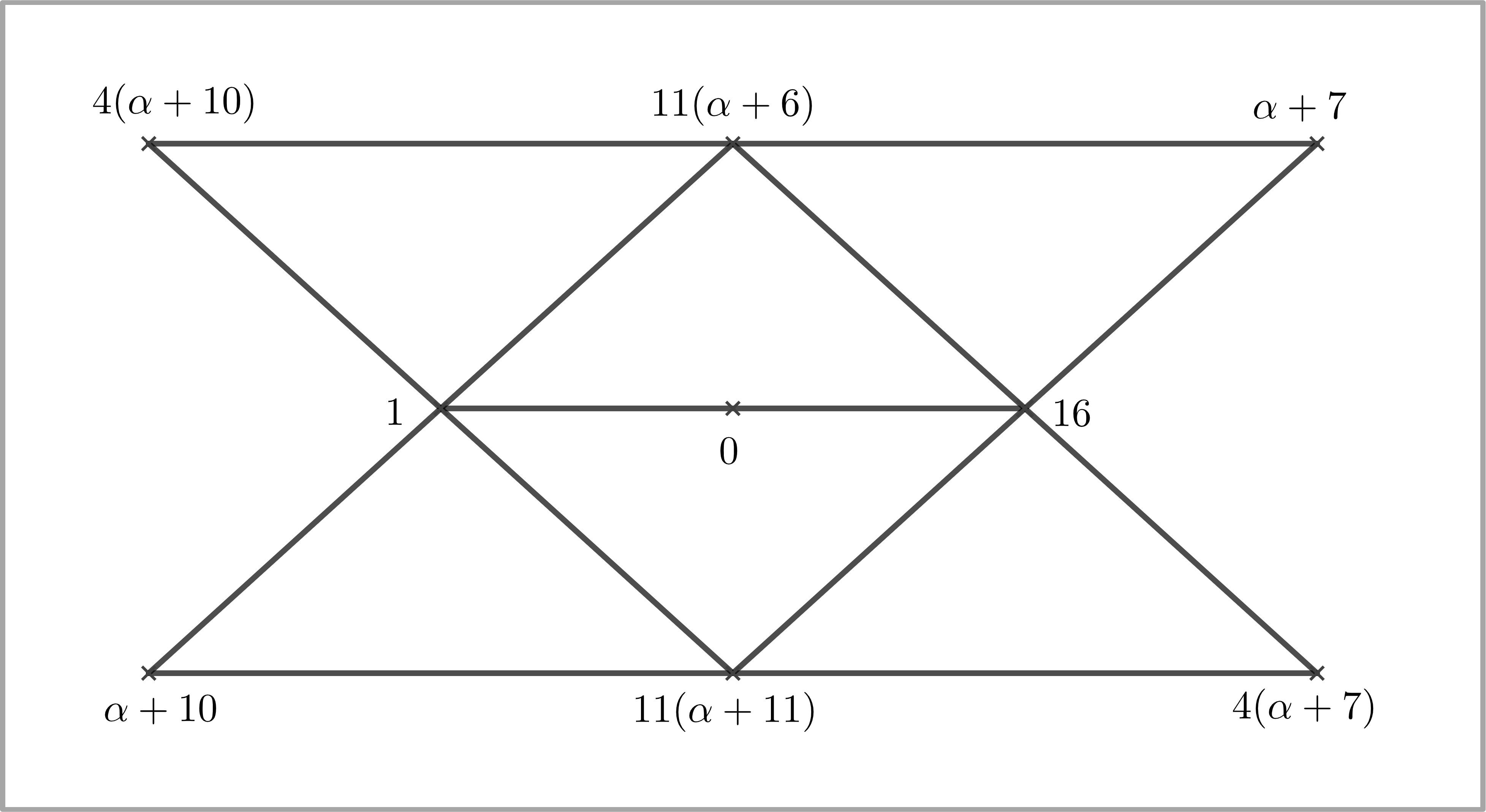}
    \quad
     \caption{\label{q=17-7}$C_{17}^{G}$}
    \end{figure}

  Set $C_{17}^{G}:=\{0\}\cup \{1,16\}\cup (\alpha+7)\{1,4\}\cup (\alpha+10)\{1,4\}\cup \{11(\alpha+11)\}\cup \{11(\alpha+6)\}$ be a subset of finite field $\FF_{17^2}$. By Magma, we have that $C_{17}^{G}$ is a maximal clique in $P(17^2)$ with size $\frac{q+1}{2}$ for $q=17$.
The structure of the clique $C_{17}^{G}$ is presented in Fig \ref{q=17-7}.

    \vskip 10mm
Remarks: There exists other lines 
that contain at least three points in $C_{17}^{G}$. They are:$\{0,\alpha+10,4(\alpha+10)\}$ and $\{0,\alpha+7,4(\alpha+7)\}$.

 \begin{lem}\label{lem:3.19}
Set $\mathcal{C}_{17}^{G}$ be the orbit of $Aut(P(17^2))$ acting on the cliques with $C_{17}^{G}\in \mathcal{C}_{17}^{G}$. Then $Aut(P(17^2))_{C_{17}^{G}}=\lg \sigma,\tau   \rg\cong \D_{6}$  where $\sigma(\gamma)=7(\alpha+6)\gamma+11(\alpha+11)$ and $\tau(\gamma)=7(\alpha+11)\gamma^{17}+11(\alpha+6)$ for any $\gamma\in \FF_{17^2}$. Moreover, $|\mathcal{C}_{17}^{F}|=13872$. 
\end{lem}  

\demo
Note that for any $\phi \in \Aut(P(17^2))$,
$\phi(\gamma)=a\gamma^{v}+b$, where $ a\in S$, $b\in \FF_{17^{2}}$, $v\in Gal(\FF_{17^{2}})$. Set $\mathcal{P}:=\{0,11(\alpha+11),11(\alpha+6)\}$ be the three intersecting points in the clique $C_{17}^{G}$, where for any $p_1,p_2\in\mathcal{P}$, there does not exist $\gamma\in\mathcal{C}_{17}^{G}/\mathcal{P} $, such that $p_1,p_2,\gamma$ belong to a line. Then $\phi(\mathcal{P})=\mathcal{P}$ for any $\phi\in Aut(P(17^{2}))_{G_{17}}$. Now we have the following six cases.

 If $\phi_1(0)=0$, $\phi_1(11(\alpha+11))=11(\alpha+11)$ and $\phi_1(11(\alpha+6))=11(\alpha+6)$, then we have that $a=1$, $b=0$ and $v=1$. It follows that $\phi_{1}(\gamma)=\gamma$ for any $\gamma\in \FF_{17^2}$.

 If $\phi_2(0)=0$, $\phi_2(11(\alpha+11))=11(\alpha+6)$ and $\phi_2(11(\alpha+6))=11(\alpha+11)$, then we have that $a=16$, $b=0$ and $|v|=2$.  It follows that $\phi_{2}(\gamma)=16\gamma^{17}$ for any $\gamma\in \FF_{17^2}$.
 
 If $\phi_3(0)=11(\alpha+11)$, $\phi_3(11(\alpha+11))=0$ and $\phi_3(11(\alpha+6))=11(\alpha+6)$, then we have that $a=10(\alpha+6)$, $b=11(\alpha+11)$ and $|v|=2$. It follows that $\phi_{3}(\gamma)=10(\alpha+6)\gamma^{17}+11(\alpha+11)$ for any $\gamma\in \FF_{17^2}$.

  If $\phi_4(0)=11(\alpha+11)$, $\phi_4(11(\alpha+11))=11(\alpha+6)$ and $\phi_4(11(\alpha+6))=0$, then we have that $a=7(\alpha+6)$, $b=11(\alpha+11)$ and $v=1$.  It follows that $\phi_{4}(\gamma)=7(\alpha+6)\gamma+11(\alpha+11)$ for any $\gamma\in \FF_{17^2}$.

  If $\phi_5(0)=11(\alpha+6)$, $\phi_5(11(\alpha+11))=11(\alpha+11)$ and $\phi_5(11(\alpha+6))=0$, then we have that $a=7(\alpha+11)$, $b=11(\alpha+6)$ and $|v|=2$.  It follows that $\phi_{5}(\gamma)=7(\alpha+11)\gamma^{17}+11(\alpha+6)$ for any $\gamma\in \FF_{17^2}$.
  
    If $\phi_6(0)=11(\alpha+6)$, $\phi_6(11(\alpha+11))=0$ and $\phi_6(11(\alpha+6))=11(\alpha+11)$, then we have that $a=10(\alpha+11)$, $b=11(\alpha+6)$ and $v=1$. It follows that $\phi_{6}(\gamma)=10(\alpha+11)\gamma+11(\alpha+6)$ for any $\gamma\in \FF_{17^2}$.

 Let $\sigma=\phi_{4}$, $\tau=\phi_{5}$, then we have that $\sigma^{3}=1$, $\tau^{2}=1$, $\tau^{-1}\sigma\tau=\sigma^{-1}$ then $\lg \sigma, \tau  \rg=\{1,\sigma,\sigma^{2},\tau,\tau\sigma, \tau\sigma^{2}\}=\{\phi_{1},\phi_{4},\phi_{6},\phi_{5}, \phi_{3},\phi_{2}\}\cong \D_6$. The action of $\sigma$ and $\tau$ acting on the points of clique $C_{17}^{G}$ is presented in the following table. 
 \begin{table}[htp]
\centering
\setlength{\tabcolsep}{0.5mm}
\small
\begin{tabular}{cccccccccc}
\toprule  
$\gamma$ & 0&1&16 & $\alpha+7$& $4(\alpha+7)$& $\alpha+10$& $4(\alpha+10)$&$11(\alpha+6)$  & $11(\alpha+11)$\\
	\midrule  
$\sigma(\gamma)$ & $11(\alpha+11)$ &$\alpha+10$ &$4(\alpha+7)$ &$16$ & $\alpha+7$& $4(\alpha+10)$& $1$& $0$&$11(\alpha+6)$\\
$\tau(\gamma)$ & $11(\alpha+6)$ &$\alpha+7$ &$4(\alpha+10)$ &$1$ & $\alpha+10$& $4(\alpha+7)$& $16$& $0$&$11(\alpha+11)$\\
\midrule 
\end{tabular}
\end{table}

Followed by the above arguments, we know that $\sigma(C_{17}^{G})=C_{17}^{G}$ and $\tau(C_{17}^{G})=C_{17}^{G}$, then $Aut(P(17^2))_{C_{17}^{G}}=\lg \sigma, \tau  \rg\cong \D_{6}$  where $\sigma(\gamma)=7(\alpha+6)\gamma+11(\alpha+11)$ and $\tau(\gamma)=7(\alpha+11)\gamma^{17}+11(\alpha+6)$ with $\gamma\in \FF_{17^2}$.

Note that $\left | Aut(P(17^{2})) \right |=\frac{q^2-1}{2}\times q^{2}\times 2=83232$, so we have that $|\mathcal{C}_{17}^{G}|=|Aut(P(17^2)):Aut(P(17^2))_{C_{17}^{G}}|=13872$.
\qed
\vskip 5mm
Set $C_{s,\gamma}:=\left \{ sx +\gamma   \mid  x\in C_{17}^{G}\cap \FF_{17^2}  \right \} $ where $s\in S$, $\gamma \in \FF_{17^{2}}$. 

\begin{lem}
Set $\mathcal{C}_{17}^{G}$ be the orbit of $Aut(P(17^2))$ acting on the cliques with $C_{17}^{G}\in \mathcal{C}_{17}$. Then  $\mathcal{C}_{17}^{G}=\{C_{s,\gamma}\mid s\in S,\gamma\in\FF_{17^2}\}$. 
\end{lem}
\demo
It is obvious that $\{C_{s,\gamma}\mid s\in S, \gamma\in\FF_{17^2}\} \subset \mathcal{C}_{17}^{G}$. If $\gamma=\gamma'=0$, now we will prove that $\C_{s,0}=C_{s',0}$ if and only if $s=s'$,  where $s,s'\in S$.

 If $C_{s,0}=C_{s',0}$, then $\{0,11(\alpha+11)s,11(\alpha+6)s\}=\{0,11(\alpha+11)s',11(\alpha+6)s'\}$. If $11(\alpha+11)s=11(\alpha+6)s'$ and $11(\alpha+6)s=11(\alpha+11)s'$ then $(\alpha+6)^2=(\alpha+11)^{2}$, a contradiction. Now, we have $s=s'$.
 
From Lemma \ref{lem:3.19}, we have that $C_{17}^{G}=\phi_4(C_{17}^{G})=\phi_6(C_{17}^{G})=7(\alpha+6)C_{17}^{G}+11(\alpha+11)=10(\alpha+11)C_{17}^{G}+11(\alpha+6)$, and $7(\alpha+6), 10(\alpha+10) \in S$.

And then $\mathcal{C}_{17}^{G}=\{C_{s,\gamma}\mid s\in S, \gamma\in\FF_{17^2}\}$, because $|\{C_{s,\gamma}\mid s\in S,\gamma\in\FF_{17^2}\}|=\frac{q^{2}-1}{2}\times q^{2}\times \frac{1}{3}=13872$.
\qed

\subsection{Maximal cliques in Paley graph $P(11^2)$}
\begin{figure}[htp]
    \centering
    \includegraphics[height=4.4cm, width=9cm]{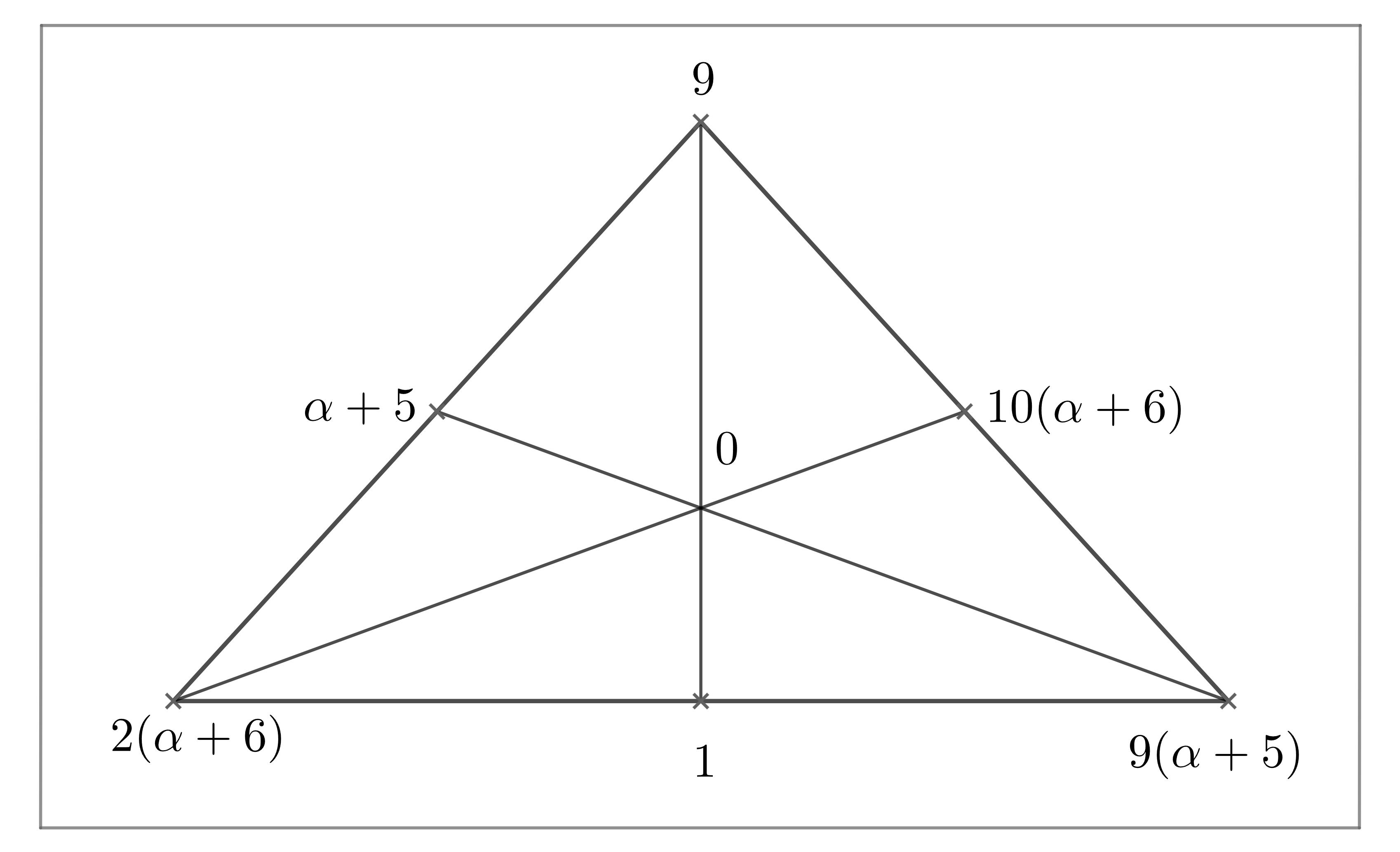}
    \quad
     \caption{\label{q=11}$C_{11}$}

    \end{figure}

We can choose a primitive element $\delta\in \FF_{11}$ and $\alpha\in\FF_{{11}^2}$ be a root of  the irreducible polynomial $x^2-\delta$ over $\FF_q$,
such that $\FF_{11}^*=\lg\delta\rg=\lg 2 \rg$ and $S_0=\{1,\alpha, \alpha+4, \alpha-4, \alpha+5, \alpha-5\}$.
Let $H=\{ 1,9 \}$  be a subset in $\FF_{q}^{\ast}$. Set $C_{11}:=\{0\}\cup H \cup (\alpha+5)H\cup 10(\alpha+6)H$ be a subset of finite field $\FF_{11^2}$. By Magma, we have that $C_{11}$ is a maximal clique in $P(11^2)$ with size $\frac{q+3}{2}$ for $q=11$.
The structure of the clique $C_{11}$ is presented in Fig \ref{q=11}.

\begin{lem}
Set $\mathcal{C}_{11}$ be the orbit of $Aut(P(11^2))$ acting on the cliques with $C_{11}\in \mathcal{C}_{11}$. Then $Aut(P(11^2))_{C_{11}}=\lg \sigma,\tau \rg\cong \D_6$, where $\sigma(\gamma)=(\alpha+5)\gamma$ and  $\tau(\gamma)=(\alpha+5)\gamma^{11}$ for any $\gamma\in \FF_{11^2}$.  Moreover,  $|\mathcal{C}_{11}|=2420$.
\end{lem}

\demo
Note that for any $\phi \in \Aut(P(11^2))$,
$\phi(\gamma)=a\gamma^{v}+b$, where $ a\in S$, $b\in \FF_{11^{2}}$, $v\in Gal(\FF_{11^{2}})$.  Set $\mathcal{P}:=\{1,\alpha+5,10(\alpha+6)\}$ be the three special points which are the intersecting points of two lines in clique $C_{11}$. Then $\phi(\mathcal{P})=\mathcal{P}$. Now we have the following six cases.

 If $\phi_1(1)=1$, $\phi_1(\alpha+5)=\alpha+5$ and $\phi_1(10(\alpha+6))=10(\alpha+6)$, then we have that $a=1$, $b=0$ and $v=1$. It follows that $\phi_{1}(\gamma)=\gamma$ for any $\gamma\in \FF_{11^2}$.

 If $\phi_2(1)=1$, $\phi_2(\alpha+5)=10(\alpha+6)$ and $\phi_2(10(\alpha+6))=\alpha+5$, then we have that $a=1$, $b=0$ and $|v|=2$. It follows that $\phi_{2}(\gamma)=\gamma^{11}$ for any $\gamma\in \FF_{11^2}$. 
 
 If $\phi_3(1)=\alpha+5$, $\phi_3(\alpha+5)=10(\alpha+6)$ and $\phi_3(10(\alpha+6))=1$, then we have that $a=\alpha+5$, $b=0$ and $v=1$. It follows that $\phi_{3}(\gamma)=(\alpha+5)\gamma$ for any $\gamma\in \FF_{11^2}$.

 If $\phi_4(1)=\alpha+5$, $\phi_4(\alpha+5)=1$ and $\phi_4(10(\alpha+6))=10(\alpha+6)$,  then we have that $a=\alpha+5$, $b=0$ and $|v|=2$. It follows that $\phi_{4}(\gamma)=(\alpha+5)\gamma^{11}$ for any $\gamma\in \FF_{11^2}$.

If $\phi_5(1)=10(\alpha+6)$, $\phi_5(\alpha+5)=\alpha+5$ and $\phi_5(10(\alpha+6))=1$,  then we have that $a=10(\alpha+6)$, $b=0$ and $|v|=2$. It follows that $\phi_{5}(\gamma)=10(\alpha+6)\gamma^{11}$ for any $\gamma\in \FF_{11^2}$.

If $\phi_6(1)=10(\alpha+6)$, $\phi_6(\alpha+5)=1$ and $\phi_6(10(\alpha+6))=\alpha+5$,  then we have that $a=10(\alpha+6)$, $b=0$ and $v=1$. It follows that $\phi_{6}(\gamma)=10(\alpha+6)\gamma$ for any $\gamma\in \FF_{11^2}$.

 Set $\sigma=\phi_{3}$, $\tau=\phi_{4}$, then we have that $\sigma^{3}=1$, $\tau^{2}=1$, $\tau^{-1}\sigma\tau=\sigma^{-1}$, then $\lg \sigma, \tau  \rg=\{1,\sigma,\sigma^{2},\tau,\tau\sigma, \tau\sigma^{2}\}=\{\phi_{1},\phi_{3},\phi_{6},\phi_{4}, \phi_{2},\phi_{5}\}\cong \D_6$. The action of $\sigma$ and $\tau$ on the points of clique $C_{11}$ is presented in the following table. 
 
 \begin{table}[htp]
\centering
\small
\begin{tabular}{cccccccc}
\toprule  
$\gamma$ & 0&1&9 & $\alpha+5$& $9(\alpha+5)$& $2(\alpha+6)$& $10(\alpha+6)$\\
	\midrule  
$\sigma(\gamma)$ & $0$ &$\alpha+5$ &$9(\alpha+5)$ &$10(\alpha+6)$ & $2(\alpha+6)$& $9$& $1$\\
$\tau(\gamma)$ & $0$ &$\alpha+5$ &$9(\alpha+5)$ &$1$ & $9$& $2(\alpha+6)$& $10(\alpha+6)$\\
\midrule 
\end{tabular}
\end{table}

Followed by the above arguments, we know that $\sigma(C_{11})=C_{11}$ and $\tau(C_{11})=C_{11}$, then $Aut(P(11^2))_{C_{11}}=\lg \sigma, \tau  \rg\cong D_{6}$  where $\sigma(\gamma)=(\alpha+5)\gamma$,  $\tau(\gamma)=(\alpha+5)\gamma^{11}$ with $\gamma\in \FF_{11^2}$.

Note that $\left | Aut(P(11^{2})) \right |=\frac{q^2-1}{2}\times q^{2}\times 2=14520$, so we have that $|\mathcal{C}_{11}|=|Aut(P(11^2)):Aut(P(11^2))_{C_{11}}|=2420$.
\qed

\vskip 5mm
Set $C_{s,\eta,\gamma}:=\{s\eta x+\gamma| x\in C_{11}\cap \FF_{11^2}\}$ for $\eta\in\{1,\alpha\}$, where $s\in \FF_{11}^{\ast}$ and $\gamma\in\FF_{11^2}$ . 

\begin{lem}
Set $\mathcal{C}_{11}$ be the orbit of $Aut(P(11^2))$ acting on the cliques with $C_{11}\in \mathcal{C}_{11}$. Then $\mathcal{C}_{11}=\{C_{s,1,\gamma} \mid s\in \FF_{11}^{\ast}, \gamma\in\FF_{11^2}\}\cup \{C_{s,\alpha,\gamma} \mid s\in \FF_{11}^{\ast}, \gamma\in\FF_{11^2}\}$.
\end{lem}
\demo 
It is obvious that $\{C_{s,1,\gamma} \mid s\in \FF_{11}^{\ast}, \gamma\in\FF_{11^2}\}\cup \{C_{s,\alpha,\gamma} \mid s\in \FF_{11}^{\ast}, \gamma\in\FF_{11^2}\}\subset \mathcal{C}_{11}$. Now we will prove that $C_{s,\eta,\gamma}=C_{s',\eta',\gamma'}$ if and only if $s=s'$, $\gamma=\gamma'$ and $\eta=\eta'$, where $s,s'\in \FF_{11}^{\ast}$, $\gamma,\gamma'\in\FF_{11^2}$ and $\eta,\eta'\in\{1,\alpha\}$.

 If $C_{s,\eta,\gamma}=C_{s',\eta',\gamma'}$, then $\gamma=\gamma'$ and the three special intersecting points $\{ s\eta+\gamma,(\alpha+5)s\eta+\gamma,10(\alpha+6)s\eta+\gamma\}=\{s' \eta'+\gamma',(\alpha+5)s'\eta'+\gamma',10(\alpha+6)s'\eta'+\gamma'\}$.  If $ s\eta+\gamma=s' \eta'+\gamma'$, then $s=s'$ and $\eta=\eta'$; if $ s\eta+\gamma=(\alpha+5)s'\eta'+\gamma'$, then we can get a contradiction; if $ s\eta+\gamma=10(\alpha+6)s'\eta'+\gamma'$, then we can get a contradiction. Now, we have $s=s'$ and $\eta=\eta'$.
 
And then $\mathcal{C}_{11}=\{C_{s,1,\gamma} \mid s\in \FF_{11}^{\ast}, \gamma\in\FF_{11^2}\}\cup \{C_{s,\alpha,\gamma} \mid s\in \FF_{11}^{\ast}, \gamma\in\FF_{11^2}\}$, because $|\{C_{s,1,\gamma} \mid s\in \FF_{11}^{\ast}, \gamma\in\FF_{11^2}\}\cup \{C_{s,\alpha,\gamma} \mid s\in \FF_{11}^{\ast}, \gamma\in\FF_{11^2}\}|=(q-1)\times q^{2}\times 2=2420$.
\qed

\subsection{Maximal cliques in Paley graph $P(19^2)$}
We can choose a primitive element $\delta\in \FF_{19}$ and $\alpha\in\FF_{{19}^2}$ be a root of  the irreducible polynomial $x^2-\delta$ over $\FF_q$,
such that $\FF_{19}^*=\lg\delta\rg=\lg 13 \rg$ and $S_0=\{1,\alpha, \alpha+2, \alpha+3, \alpha+4, \alpha+6, \alpha-2, \alpha-3, \alpha-4, \alpha-6\}$.
\subsubsection{$C_{19}^{A}$-construction}
\begin{figure}[htp]
    \centering
    \includegraphics[height=4.4cm, width=9cm]{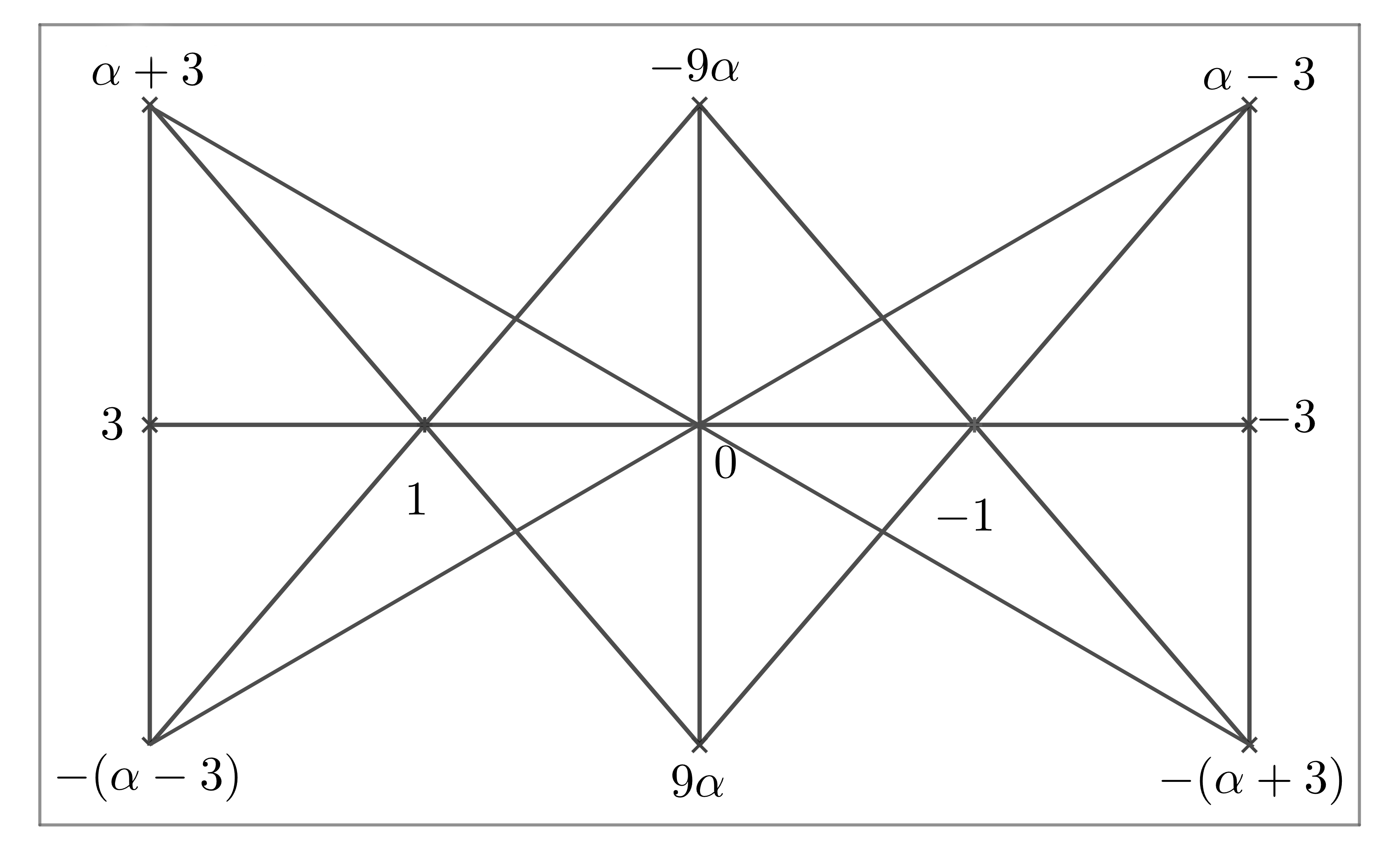}
    \quad
     \caption{\label{q=19-1}$C_{19}^{A}$}
     Remarks: There exists other lines 
that contain at least three points in $C_{19}^{A}$. 
     
     They are:$\{3,-9\alpha,\alpha-3\}$,$\{3,9\alpha,-(\alpha+3)\}$,$\{-3,-9\alpha,\alpha+3\}$,$\{-3,9\alpha,-(\alpha-3)\}$.
    \end{figure}

Let $H=\{1,-1 \}$ be a subset in $\FF_{19}^{\ast}$. Set $C_{19}^{A}:= \{ 0  \}\cup   \{H,3H \} \cup 9\alpha H  \cup  ( \alpha +3 ) H \cup ( \alpha -3 ) H$ be a subset of finite field $\FF_{{19}^2}$. By Magma, we have that $C_{19}^{A}$ is a maximal clique in $P(19^2)$ with size $\frac{q+3}{2}$ for $q=19$.
The structure of the clique $C_{19}^{A}$ is presented in Fig \ref{q=19-1}.

\begin{lem}
 Set $\mathcal{C}_{19}^{A}$ be the orbit of $Aut(P(19^2))$ acting on the cliques with $C_{19}^{A}\in \mathcal{C}_{19}^{A}$. Then $Aut(P(19^2))_{C_{19}^{A}}=\lg \sigma\rg\times\lg \tau\rg\cong \ZZ_{2}\times\ZZ_{2}$, where $\sigma(\gamma)=-\gamma$,  $\tau(\gamma)=\gamma^{19}$ for any $\gamma\in \FF_{19^2}$. Moreover, $|\mathcal{C}_{19}^{A}|=32490$.
\end{lem}
\demo
 Note that for any $\phi \in \Aut(P(19^2))$,
$\phi(\gamma)=a\gamma^{v}+b$, where $ a\in S$, $b\in \FF_{19^{2}}$, $v\in Gal(\FF_{19^{2}})$. Set $\mathcal{P}:=\{0,3,-3\}$ be the three special points which are the intersecting points of four lines in Figure \ref{q=19-1}, and one of these lines contain five points in the clique.  Then $\phi(\mathcal{P})=\mathcal{P}$. Now we have the following six cases.

 If $\phi(0)=0$, $\phi(3)=3$ and $\phi(-3)=-3$, then we have that $a=1$, $b=0$ and $ v\in Gal(\FF_{19^{2}})$. Let $\phi_1(\gamma)=\gamma$ and $\phi_2(\gamma)=\gamma^{19}$ for any for any $\gamma\in \FF_{19^2}$.

 If $\phi(0)=0$, $\phi(3)=-3$ and $\phi(-3)=3$, then we have that $a=-1$, $b=0$ and $ v\in Gal(\FF_{19^{2}})$.  Let $\phi_3(\gamma)=-\gamma$ and $\phi_4(\gamma)=-\gamma^{19}$ for any for any $\gamma\in \FF_{19^2}$.

  If $\phi(0)=3$, $\phi(3)=0$ and $\phi(-3)=-3$, then there are no solutions for $a,b$ and $v$.
 If $\phi(0)=3$, $\phi(3)=-3$ and $\phi(-3)=0$, then there are no solutions for $a,b$ and $v$. 
 If $\phi(0)=-3$, $\phi(3)=0$ and $\phi(-3)=3$, then there are no solutions for $a,b$ and $v$.
 If $\phi(0)=-3$, $\phi(3)=3$ and $\phi(-3)=0$, then there are no solutions for $a,b$ and $v$.
  
Set $\sigma=\phi_{3}$, $\tau=\phi_{2}$. Then we have that $\sigma^{2}=1$, $\tau^{2}=1$ and $\lg \sigma\rg\times\lg \tau\rg=\{1,\sigma,\tau,\sigma\tau\}=\{\phi_{1},\phi_{3},\phi_{2},\phi_{4}\}\cong \ZZ_{2}\times\ZZ_{2}$. The action of $\sigma$ and $\tau$ on the points in the clique $C_{19}^{A}$ is presented in the following table.
  
 \begin{table}[htp]
\centering

\begin{tabular}{cccccccccccc}
\toprule  
$\gamma$ & 0&1&-1 & $3$& $-3$& $9\alpha$& $-9\alpha$&$(\alpha+3)$  & $-(\alpha+3)$&$\alpha-3$&$-(\alpha-3)$\\
	\midrule  
$\sigma(\gamma)$ & 0&-1&1 & $-3$& $3$& $-9\alpha$& $9\alpha$&$-(\alpha+3)$ & $(\alpha+3)$&$-(\alpha-3)$&$(\alpha-3)$\\
$\tau(\gamma)$ & 0&1&-1 & $3$& $-3$& $-9\alpha$& $9\alpha$&$-(\alpha-3)$ & $(\alpha-3)$&$-(\alpha+3)$&$(\alpha+3)$\\
\midrule 
\end{tabular}
\end{table}

Followed by the above arguments, we know that $\sigma(C_{19}^{A})=C_{19}^{A}$ and $\tau(C_{19}^{A})=C_{19}^{A}$, then $Aut(P(19^2))_{C_{19}^{A}}=\lg \sigma\rg\times\lg \tau\rg\cong \ZZ_{2}\times\ZZ_{2}$ where $\sigma(\gamma)=-\gamma$,  $\tau(\gamma)=\gamma^{19}$ for any $\gamma\in \FF_{19^2}$.

Note that $\left | Aut(P(19^{2})) \right |=\frac{q^2-1}{2}\times q^{2}\times 2=129960$, so we have that $|\mathcal{C}_{19}^{A}|=|Aut(P(19^2)):Aut(P(19^2))_{C_{19}^{A}}|=32490$.
\qed

\vskip 5mm
Set $C_{s,\gamma,i}:=\left \{ six +\gamma   \mid  x\in C_{19}^{A}\cap \FF_{19^2}  \right \} $ where $s\in S_{0}$, $i\in \{ 1\dots9\}$, $\gamma \in \FF_{19^{2}}$. 

\begin{lem}
Set $\mathcal{C}_{19}^{A}$ be the orbit of $Aut(P(19^2))$ acting on the cliques with $C_{19}^{A}\in \mathcal{C}_{19}$. Then  $\mathcal{C}_{19}^{A}=\{C_{s,\gamma,i}\mid s\in S_{0},i\in \{1\dots9\}, \gamma\in\FF_{19^2}\}$.
\end{lem}
\demo
It is obvious that $\{C_{s,\gamma,i}\mid s\in S_{0},i\in \{1\dots9\}, \gamma\in\FF_{19^2}\} \subset \mathcal{C}_{19}^{A}$. Now we will prove that $\C_{s,\gamma,i}=C_{s',\gamma',i'}$ if and only if $s=s'$, $\gamma=\gamma'$ and $i=i'$, where $s,s'\in S_{0}$, $i,i'\in\{1\dots9\}$ and $\gamma,\gamma'\in\FF_{19^2}$.

If $C_{s,\gamma,i}=C_{s',\gamma',i'}$, then $\gamma=\gamma'$ and the subset of a line $siH'+\gamma=s'i'H'+\gamma'$. Then $si\in\{s'i',-s'i'\}$, and we have that $s=s'$ and $i=i'$.

    And then $\mathcal{C}_{19}^{A}=\{C_{s,\gamma,i}\mid s\in S_{0}, i\in \{1\dots9\}, \gamma\in\FF_{19^2}\}$, because $|\{C_{s,\gamma,i}\mid s\in S_{0},i\in \{1\dots9\}, \gamma\in\FF_{19^2}\}|=\frac{q+1}{2}\times q^{2}\times 9=32490$.
\qed

\subsubsection{$C_{19}^{B}$-construction}

\begin{figure}[htp]
    \centering
    \includegraphics[height=4.4cm, width=9cm]{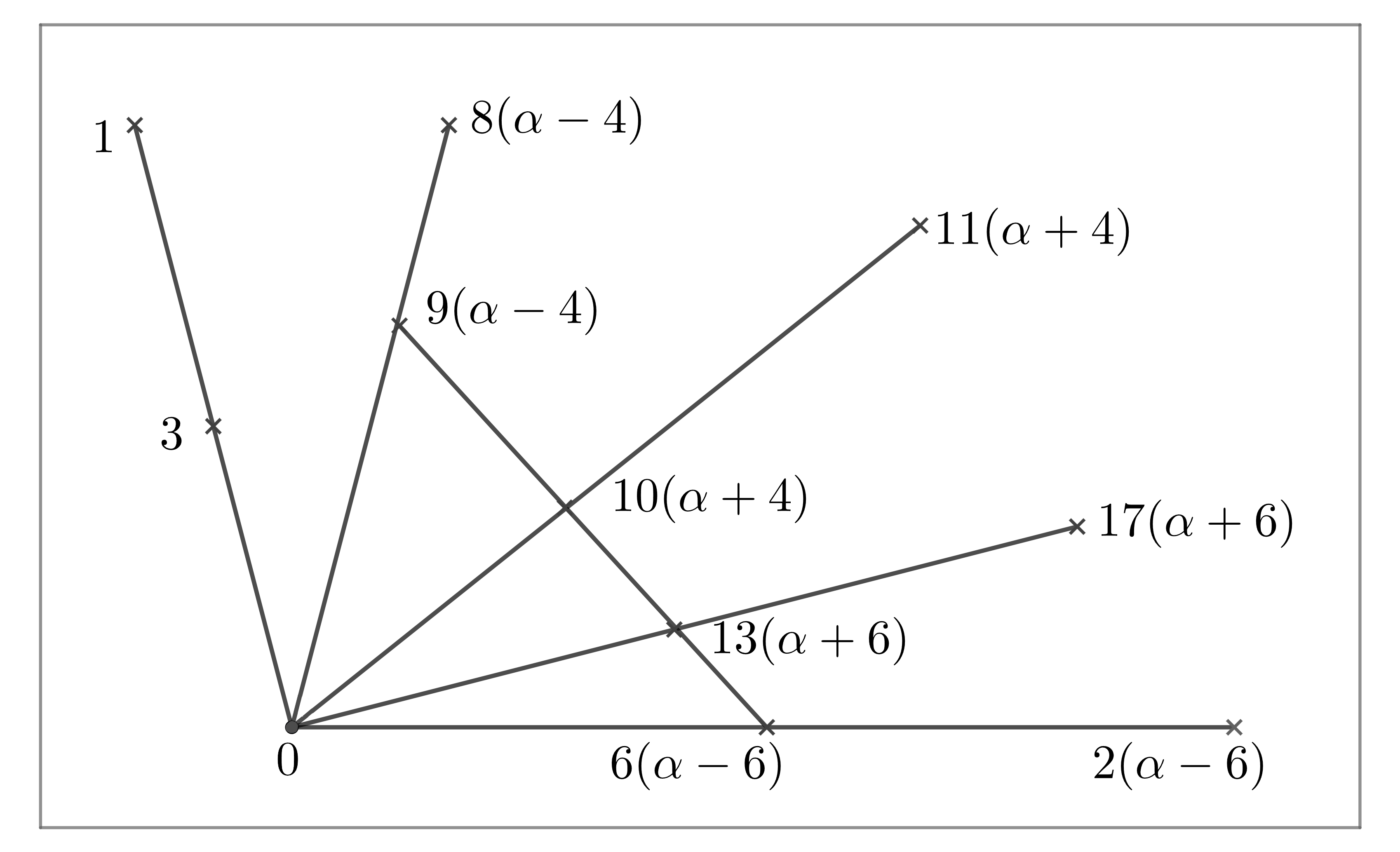}
    \quad
     \caption{\label{q=19-2}$C_{19}^{B}$}
    \end{figure}

Let $H=\left \{1,3 \right \}$ be a subset $\FF_{19}^{\ast}$, 
Set $C_{19}^{B}:= \{ 0 \}\cup  H \cup 2H ( \alpha -6 ) \cup 17H ( \alpha +6  )\cup 9H ( \alpha -4  ) \cup10 H( \alpha +4  )$ be a subset of finite field $\FF_{19^2}$. By Magma, we have that $C_{19}^{B}$ is a maximal clique in $P(19^2)$ with size $\frac{q+3}{2}$ for $q=19$.
The structure of the clique $C_{19}^{B}$ is presented in Fig \ref{q=19-2}.

Remarks: There exists other lines 
that contain at least three points in $C_{19}^{B}$.  
     They are:$\{1,13(\alpha+6),2(\alpha-6),8(\alpha-4)\}$,$\{1,6(\alpha-6),17(\alpha+6),11(\alpha+4)\}$,$\{3,11(\alpha+4),9(\alpha-4),2(\alpha-6)\}$,$\{3,8(\alpha-4),10(\alpha+4),17(\alpha+6)\}$.
\begin{lem}\label{lem:3.25}
 Set $\mathcal{C}_{19}^{B}$ be the orbit of $Aut(P(19^2))$ acting on the cliques with $C_{19}^{B}\in \mathcal{C}_{19}^{B}$. Then $Aut(P(19^2))_{C_{19}^{B}}=\lg \sigma, \tau  \rg\cong \D_{10}$, where $ \sigma(\gamma)=10(\alpha+4)\gamma$ and $\tau(\gamma)=10(\alpha+4)\gamma^{19}$ for any  $\gamma\in \FF_{{19}^2}$.  Moreover, $|\mathcal{C}_{19}^{B}|=12996$.    
\end{lem}
\demo
Note that for any $\phi \in \Aut(P(19^2))$,
$\phi(\gamma)=a\gamma^{v}+b$, where $ a\in S$, $b\in \FF_{19^{2}}$, $v\in Gal(\FF_{19^{2}})$. 
Note that the point $0$ in the clique $C_{19}^B$ is a intersecting point of five lines, which contain three points in this clique.
If $\phi\in Aut(P(19^{2}))_{C_{19}^B}$, then $\phi(0)=0 $ and $b=0$. Set $\mathcal{L}:=\{H,2(\alpha-6)H,17(\alpha+6)H,9(\alpha-4)H,10(\alpha+4)H\}$ be the subset of five lines in clique $C_{19}^{B}$,  which are presented in Fig \ref{q=19-2}. Then $\phi(\mathcal{L})=\mathcal{L}$. 

 If $\phi(H)=H$, then we have that $a=1$, $b=0$. It follows that $\phi(\gamma)=\gamma^v$. Set $\phi_{1}(\gamma)=\gamma$ and $\phi_{2}(\gamma)=\gamma^{19}$ for any $\gamma\in \FF_{q^2}$.

 If $\phi(H)=2(\alpha-6)H$, then we have that $a=2(\alpha-6)$, $b=0$. It follows that  $\phi(\gamma)=2(\alpha-6)\gamma^v$. Set $\phi_{3}(\gamma)=2(\alpha-6)\gamma$ and $\phi_{4}(\gamma)=2(\alpha-6)\gamma^{19}$ for any $\gamma\in \FF_{q^2}$.  

 Similarly, we have that $\phi_{5}(\gamma)=17(\alpha+6)\gamma$, $\phi_{6}(\gamma)=17(\alpha+6)\gamma^{19}$, $\phi_{7}(\gamma)=9(\alpha-4)\gamma$, $\phi_{8}(\gamma)=9(\alpha-4)\gamma^{19}$
 , $\phi_{9}(\gamma)=10(\alpha+4)\gamma$, $\phi_{10}(\gamma)=10(\alpha+4)\gamma^{19}$.

 Set $\sigma=\phi_{9}$, $\tau=\phi_{10}$, then we have that $\sigma^{5}=1$, $\tau^{2}=1$, $\tau^{-1}\sigma\tau=\sigma^{-1}$, then $\lg \sigma, \tau  \rg=\{1,\sigma,\sigma^2,\sigma^3,\sigma^4,\tau,\tau\sigma,\tau\sigma^2,\tau\sigma^3,\tau\sigma^4 \}=\{\phi_{1},\phi_{9},\phi_{3},\phi_{5},\phi_{7},\phi_{10},\phi_{2},\phi_{8}, \phi_{6},\phi_{4},\}\cong \D_{10}$. The action of $\sigma$ and $\tau$ acting on the points in clique $C_{19}^{B}$ is presented by the following table.

 \begin{table}[htp]
\centering
\setlength{\tabcolsep}{0.6mm}
\begin{tabular}{cccccccccccc}
\toprule  

$\gamma$ & 0&1&3 & $6(\alpha-6)$& $2(\alpha-6)$& $13(\alpha+6)$& $17(\alpha+6)$  & $10(\alpha+4)$ &$11(\alpha+4)$&$9(\alpha-4)$&$8(\alpha-4)$\\
	\midrule  
$\sigma(\gamma)$ &0&$10(\alpha+4)$&$11(\alpha+4)$&$13(\alpha+6)$ & $17(\alpha+6)$& $8(\alpha-4)$& $9(\alpha-4)$& $2(\alpha-6)$&$6(\alpha-6)$ & $1$&$3$\\
$\tau(\gamma)$ &0&$10(\alpha+4)$&$11(\alpha+4)$&$8(\alpha-4)$ & $9(\alpha-4)$& $13(\alpha+6)$& $17(\alpha+6)$& $1$&$3$& $2(\alpha-6)$&$6(\alpha-6)$ \\

\midrule 
\end{tabular}
\end{table}

Followed by the above arguments, we know that $\sigma(C_{19}^{B})=C_{19}^{B}$ and $\tau(C_{19}^{B})=C_{19}^{B}$. Then $Aut(P(19^2))_{C_{19}^{B}}=\lg \sigma, \tau  \rg\cong \D_{10}$,  where $ \sigma(\gamma)=10(\alpha+4)\gamma$ and $\tau(\gamma)=10(\alpha+4)\gamma^{19}$ with $\gamma\in \FF_{19^2}$.

Note that $\left | Aut(P(19^{2})) \right |=\frac{q^2-1}{2}\times q^{2}\times 2=129960$, so we have that $|\mathcal{C}_{19}^{B}|=|Aut(P(19^2)):Aut(P(19^2))_{C_{19}^{B}}|=12996$.
\qed

\vskip 5mm
Set $C_{s,\eta,\gamma}:=\{s\eta x+\gamma| x\in C_{19}^B\cap \FF_{19^2}\}$ for $\eta\in\{1,\alpha\}$, where $s\in \FF_{19}^{\ast}$ and $\gamma\in\FF_{19^2}$ .

\begin{lem}
Set $\mathcal{C}_{19}^{B}$ be the orbit of $Aut(P(19^2))$ acting on the cliques with $C_{19}^{B}\in \mathcal{C}_{19}^{B}$. Then $\mathcal{C}_{19}^{B}=\{C_{s,1,\gamma} \mid s\in \FF_{19}^{\ast}, \gamma\in\FF_{19^2}\}\cup \{C_{s,\alpha,\gamma} \mid s\in \FF_{19}^{\ast}, \gamma\in\FF_{19^2}\}$.  
\end{lem}
\demo 
It is obvious that $\{C_{s,1,\gamma} \mid s\in \FF_{19}^{\ast}, \gamma\in\FF_{19^2}\}\cup \{C_{s,\alpha,\gamma} \mid s\in \FF_{19}^{\ast}, \gamma\in\FF_{19^2}\}\subset \mathcal{C}_{19}^B$. Now we will prove that $C_{s,\eta,\gamma}=C_{s',\eta',\gamma'}$ if and only if $s=s'$, $\gamma=\gamma'$ and $\eta=\eta'$, where $s,s'\in \FF_{19}^{\ast}$, $\gamma,\gamma'\in\FF_{19^2}$ and $\eta,\eta'\in\{1,\alpha\}$.

Note that the intersecting point of five lines in the clique $C_{s,\eta,\gamma}$ is $\gamma$. 
 If $C_{s,\eta,\gamma}=C_{s',\eta',\gamma'}$, then $\gamma=\gamma'$ and the subset of five lines $\{s\eta 
 H,2s\eta(\alpha-6) H,17s\eta(\alpha+6) H, 9s\eta(\alpha-4) H,10s\eta(\alpha+4) H\}= \{s'\eta' H,2s'\eta'(\alpha-6) H, 17s'\eta'(\alpha+6) H, 9s'\eta'(\alpha-4) H,10s'\eta'(\alpha+4) H\}$. Then $s\eta\in\{s'\eta' ,2s'\eta'(\alpha-6), 17s'\eta'(\alpha+6), 9s'\eta'(\alpha-4),10s'\eta'(\alpha+4)\}$, we have that $s=s'$ and $\eta=\eta'$.

And then $\mathcal{C}_{19}^B=\{C_{s,1,\gamma} \mid s\in \FF_{19}^{\ast}, \gamma\in\FF_{19^2}\}\cup \{C_{s,\alpha,\gamma} \mid s\in \FF_{19}^{\ast}, \gamma\in\FF_{19^2}\}$, because $|\{C_{s,1,\gamma} \mid s\in \FF_{19}^{\ast}, \gamma\in\FF_{19^2}\}\cup \{C_{s,\alpha,\gamma} \mid s\in \FF_{19}^{\ast}, \gamma\in\FF_{19^2}\}|=(q-1)\times q^{2}\times 2=12996$.
\qed

\subsection{Maximal cliques in Paley graph $P(23^2)$}
We can choose a primitive element $\delta\in \FF_{23}$ and $\alpha\in\FF_{{23}^2}$ be a root of  the irreducible polynomial $x^2-\delta$ over $\FF_q$,
such that $\FF_{23}^*=\lg\delta\rg=\lg 14 \rg$ and $S_0=\{1,\alpha, \alpha+1, \alpha+2, \alpha+6, \alpha+9, \alpha+11,\alpha-1, \alpha-2, \alpha-6, \alpha-9, \alpha-11\}$.

\subsubsection{$C_{23}^{A}$-construction}
\begin{figure}[htp]
    \centering
    \includegraphics[height=4.4cm, width=9cm]{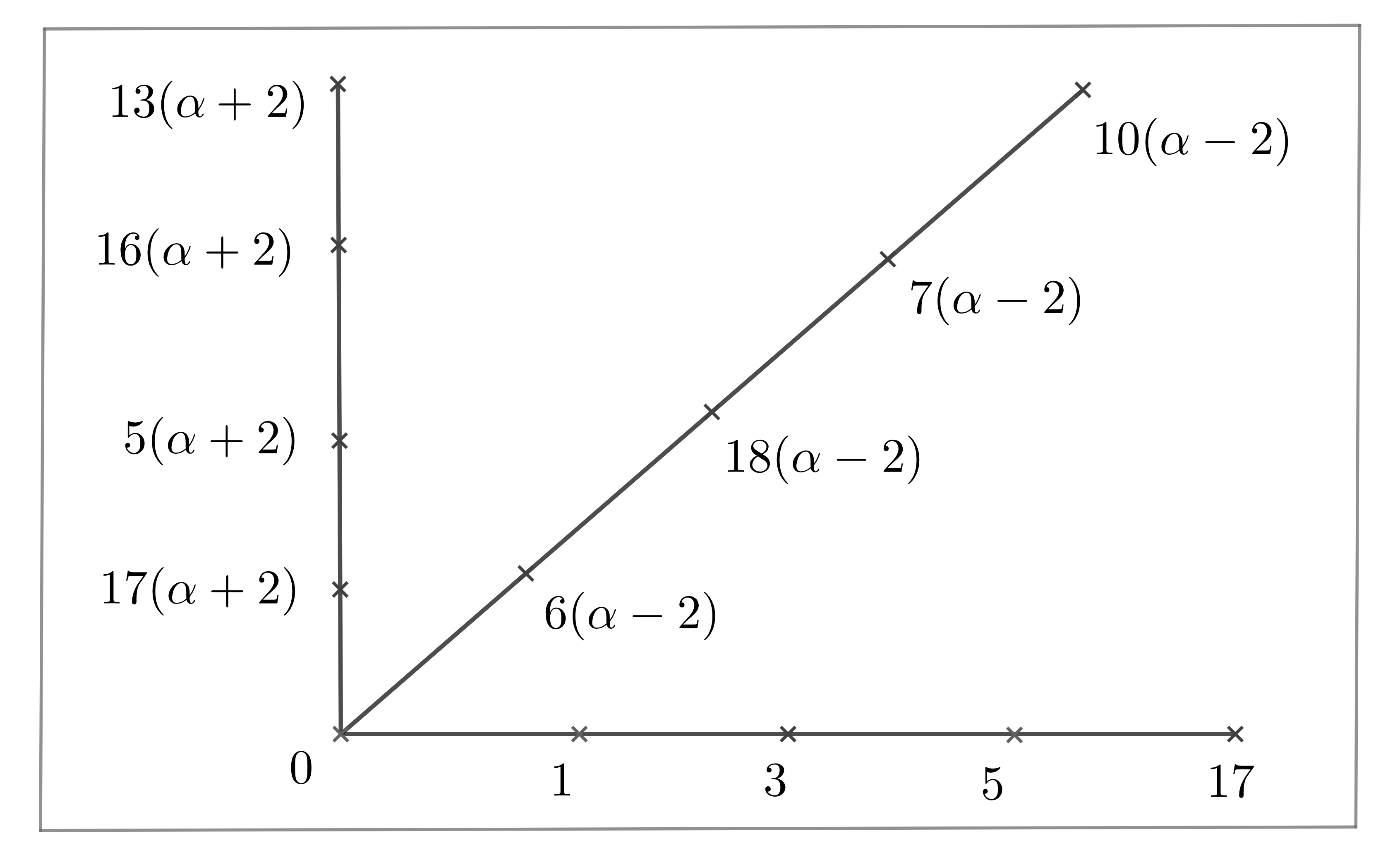}
    \quad
     \caption{\label{q=23-1}$C_{23}^{A}$}
    \end{figure}

Let H=$\left \{1,3,5,17\right \} $ be a subset $\FF_{23}^{\ast}$. Set $C_{23}^{A}:=\{ 0 \} \cup H\cup 17H ( \alpha +2  ) \cup 6H (\alpha -2 )$ be a subset of finite field $\FF_{23^2}$. By Magma, we have that $C_{23}^{A}$ is a maximal clique in $P(23^2)$ with size $\frac{q+3}{2}$ for $q=23$.
The structure of the clique $C_{23}^{A}$ is presented in Fig \ref{q=23-1}

    \vskip 10mm
  Remarks: There exists other lines 
that  contain at least three points in $C_{23}^{A}$.  
     They are:$\{1,18(\alpha-2),16(\alpha+2)\}$, $\{1,7(\alpha-2),5(\alpha+2)\}$, $\{3,6(\alpha-2),16(\alpha+2)\}$, $\{3,7(\alpha-2),17(\alpha+2)\}$, $\{3,10(\alpha-2),13(\alpha+2)\}$, $\{5,6(\alpha-2),5(\alpha+2)\}$,$\{5,18(\alpha-2),17(\alpha+2)\}$,$\{17,18(\alpha-2),13(\alpha+2)\}$,$\{17,10(\alpha-2),5(\alpha+2)\}$. 
\begin{lem}
Set $\mathcal{C}_{23}^{A}$ be the orbit of $Aut(P(23^2))$ acting on the cliques with $C_{23}^{A}\in \mathcal{C}_{23}^{A}$. Then $Aut(P(23^2))_{C_{23}^{A}}=\lg \sigma, \tau\rg\cong  \D_{6}$, where $\sigma(\gamma)=6(\alpha+21)\gamma$ and $\tau(\gamma)=6(\alpha+21)\gamma^{23}$ for any  $\gamma\in \FF_{23^2}$.  Moreover, $|\mathcal{C}_{23}^{A}|=46552$.
\end{lem}
\demo
Note that for any $\phi \in \Aut(P(23^2))$,
$\phi(\gamma)=a\gamma^{v}+b$, where $ a\in S$, $b\in \FF_{23^{2}}$, $v\in Gal(\FF_{23^{2}})$. 
Note that the point $0$ in the clique $C_{23}^A$ is a intersecting point of three lines, which contain five points in this clique.
If $\phi\in Aut(P(23^{2}))_{C_{23}^{A}}$, then  $\phi(0)=0$ and $b=0$. Set $\mathcal{L}:=\{H,17(\alpha+2)H,6(\alpha-2)H\}$ be the subset of  three lines, which are presented in Fig \ref{q=23-1}. Then $\phi(\mathcal{L})=\mathcal{L}$. Now we have the following six cases.

If $\phi_1(H)=H$, $\phi_1(17(\alpha+2)H)=17(\alpha+2)H$ and $\phi_1(6(\alpha-2)H)=6(\alpha-2)H$, then we have that $a=1$ and $v=1$.  It follows that $\phi_{1}(\gamma)=\gamma$ for any $\gamma\in \FF_{23^2}$.

If  $\phi_2(H)=H$, $\phi_2(17(\alpha+2)H)=6(\alpha-2)H$ and $\phi_2(6(\alpha-2)H)=17(\alpha+2)H$, then we have that $a=1$ and $|v|=2$.  It follows that $\phi_{2}(\gamma)=\gamma^{23}$ for any $\gamma\in \FF_{23^2}$. 

If  $\phi_3(H)=17(\alpha+2)H$, $\phi_3(17(\alpha+2)H)=6(\alpha-2)H$ and $\phi_3(6(\alpha-2)H)=H$, then we have that $a=17(\alpha+2)$ and $v=1$.  It follows that $\phi_{3}(\gamma)=17(\alpha+2)\gamma$ for any $\gamma\in \FF_{23^2}$. 

If  $\phi_4(H)=17(\alpha+2)H$, $\phi_4(17(\alpha+2)H)=H$ and $\phi_4(6(\alpha-2)H)=6(\alpha-2)H$, then we have that $a=17(\alpha+2)$ and $|v|=2$.  It follows that $\phi_{4}(\gamma)=17(\alpha+2)\gamma^{23}$ for any $\gamma\in \FF_{23^2}$.

If  $\phi_5(H)=6(\alpha-2)H$, $\phi_5(17(\alpha+2)H)=H$ and $\phi_5(6(\alpha-2)H)=17(\alpha+2)H$, then we have that $a=6(\alpha-2)$ and $v=1$.  It follows that $\phi_{5}(\gamma)=6(\alpha-2)\gamma$ for any $\gamma\in \FF_{23^2}$.

If  $\phi_6(H)=6(\alpha-2)H$, $\phi_6(17(\alpha+2)H)=17(\alpha+2)H$ and $\phi_6(6(\alpha-2)H)=H$, then we have that $a=6(\alpha-2)$ and $v=23$.  It follows that $\phi_{6}(\gamma)=6(\alpha-2)\gamma^{23}$ for any $\gamma\in \FF_{23^2}$.

Set $\sigma=\phi_{5}$, $\tau=\phi_{6}$ with $\gamma\in \FF_{23^2}$, then we have that $\sigma^{3}=1$, $\tau^{2}=1$ $\tau^{-1}\sigma\tau=\sigma^{-1}$, then $\lg \sigma, \tau \rg=\{1,\sigma,\sigma^{2},\tau,\tau\sigma, \tau\sigma^{2}\}=\{\phi_{1},\phi_{5},\phi_{3},\phi_{6}, \phi_{2},\phi_{4}\}$. The action of $\sigma$ and $\tau$ acting on the points in clique $C_{23}^{A}$ is presented in the following table.

 \begin{table}[htp]
\centering
\small
\setlength{\tabcolsep}{0.45mm}
\begin{tabular}{cccccccccccccc}
\toprule  
$\gamma$ & 0&1&3&5 & $17$& $17(\alpha+2)$& $5(\alpha+2)$& $16(\alpha+2)$  & $13(\alpha+2)$ & $6(\alpha-2)$& $18(\alpha-2)$& $7(\alpha-2)$  & $10(\alpha-2)$\\
	\midrule  
$\sigma(\gamma)$ &  0&$6(\alpha-2)$& $18(\alpha-2)$& $7(\alpha-2)$  & $10(\alpha-2)$& $1$& $3$  & $5$ & $17$& $17(\alpha+2)$& $5(\alpha+2)$  & $16(\alpha+2)$&$13(\alpha
+2)$ \\
$\tau(\gamma)$&0& $6(\alpha-2)$& $18(\alpha-2)$& $7(\alpha-2)$  & $10(\alpha-2)$&  $17(\alpha+2)$& $5(\alpha+2)$& $16(\alpha+2)$  & $13(\alpha+2)$& 1& 3  & 5&17 \\

\midrule 
\end{tabular}
\end{table}

Followed by the above arguments, we know that $\sigma(C_{23}^{A})=C_{23}^{A}$ and $\tau(C_{23}^{A})=C_{23}^{A}$, then $Aut(P(23^2))_{C_{23}^{A}}=\lg \sigma,\tau\rg\cong \D_{6}$ where $\sigma(\gamma)=6(\alpha+21)\gamma$, $\tau(\gamma)=6(\alpha+21)\gamma^{23}$ with $\gamma\in \FF_{23^2}$.

Note that $\left | Aut(P(23^{2})) \right |=\frac{q^2-1}{2}\times q^{2}\times 2=279312$, so we have that $|\mathcal{C}_{23}^{A}|=|Aut(P(23^2)):Aut(P(23^2))_{C_{23}^{A}}|=46552$.
\qed

\vskip 5mm
Set $C_{s,\eta,\gamma}:=\{s\eta x+\gamma| x\in C_{23}^{A}\cap \FF_{23^2}\}$ for $\eta\in\{1,\alpha,\alpha+1,\alpha-1\}$, where $s\in \FF_{23}^{\ast}$ and $\gamma\in\FF_{23^2}$ .

\begin{lem}
Set $\mathcal{C}_{23}^{A}$ be the orbit of $Aut(P(23^2))$ acting on the cliques with $C_{23}^{A}\in \mathcal{C}_{23}^{A}$. Then $\mathcal{C}_{23}^{A}=\bigcup_{s,\eta,\gamma}\{C_{s,\eta,\gamma } \mid s\in \FF_{23}^{\ast}, \gamma\in\FF_{23^2}, \eta\in\{1,\alpha,\alpha+1,\alpha-1\} \}$.  
\end{lem}
\demo 
It is obvious that $\{C_{s,1,\gamma} \mid s\in \FF_{23}^{\ast}, \gamma\in\FF_{23^2}\}\cup \{C_{s,\alpha,\gamma} \mid s\in \FF_{23}^{\ast}, \gamma\in\FF_{23^2}\}\cup \{C_{s,\alpha+1,\gamma} \mid s\in \FF_{23}^{\ast}, \gamma\in\FF_{23^2}\}\cup \{C_{s,\alpha-1,\gamma} \mid s\in \FF_{23}^{\ast}, \gamma\in\FF_{23^2}\}\subset \mathcal{C}_{23}^A$. Now we will prove that $C_{s,\eta,\gamma}=C_{s',\eta',\gamma'}$ if and only if $s=s'$, $\gamma=\gamma'$ and $\eta=\eta'$, where $s,s'\in \FF_{23}^{\ast}$, $\gamma,\gamma'\in\FF_{{23}^2}$ and $\eta,\eta'\in\{1,\alpha,\alpha+1,\alpha-1\}$.

Note that the intersection point of three lines in the clique $C_{s,\eta,\gamma}$ is $\gamma$.
 If $C_{s,\eta,\gamma}=C_{s',\eta',\gamma'}$, then $\gamma=\gamma'$ and the subset of three lines $\{s\eta H,17s\eta(\alpha+2) H,6s\eta(\alpha-2) H\}= \{s'\eta' H,17s'\eta'(\alpha+2) H,6s'\eta'(\alpha-2) H\}$. Then $s\eta\in\{s'\eta',17s'\eta'(\alpha+2),6s'\eta'(\alpha-2)\}$, and we have that $s=s'$ and $\eta=\eta'$.

And then $\mathcal{C}_{23}^A=\bigcup_{s,\eta,\gamma}\{C_{s,\eta,\gamma } \mid s\in \FF_{23}^{\ast}, \gamma\in\FF_{23^2}, \eta\in\{1,\alpha,\alpha+1,\alpha-1\} \}$, because $|\bigcup_{s,\eta,\gamma}\{C_{s,\eta,\gamma } \mid s\in \FF_{23}^{\ast}, \gamma\in\FF_{23^2}, \eta\in\{1,\alpha,\alpha+1,\alpha-1\} \}|=(q-1)\times q^{2}\times 4=46552$.
\qed

\subsubsection{$C_{23}^{B}$-construction}

\begin{figure}[htp]
    \centering
    \includegraphics[height=4.4cm, width=9cm]{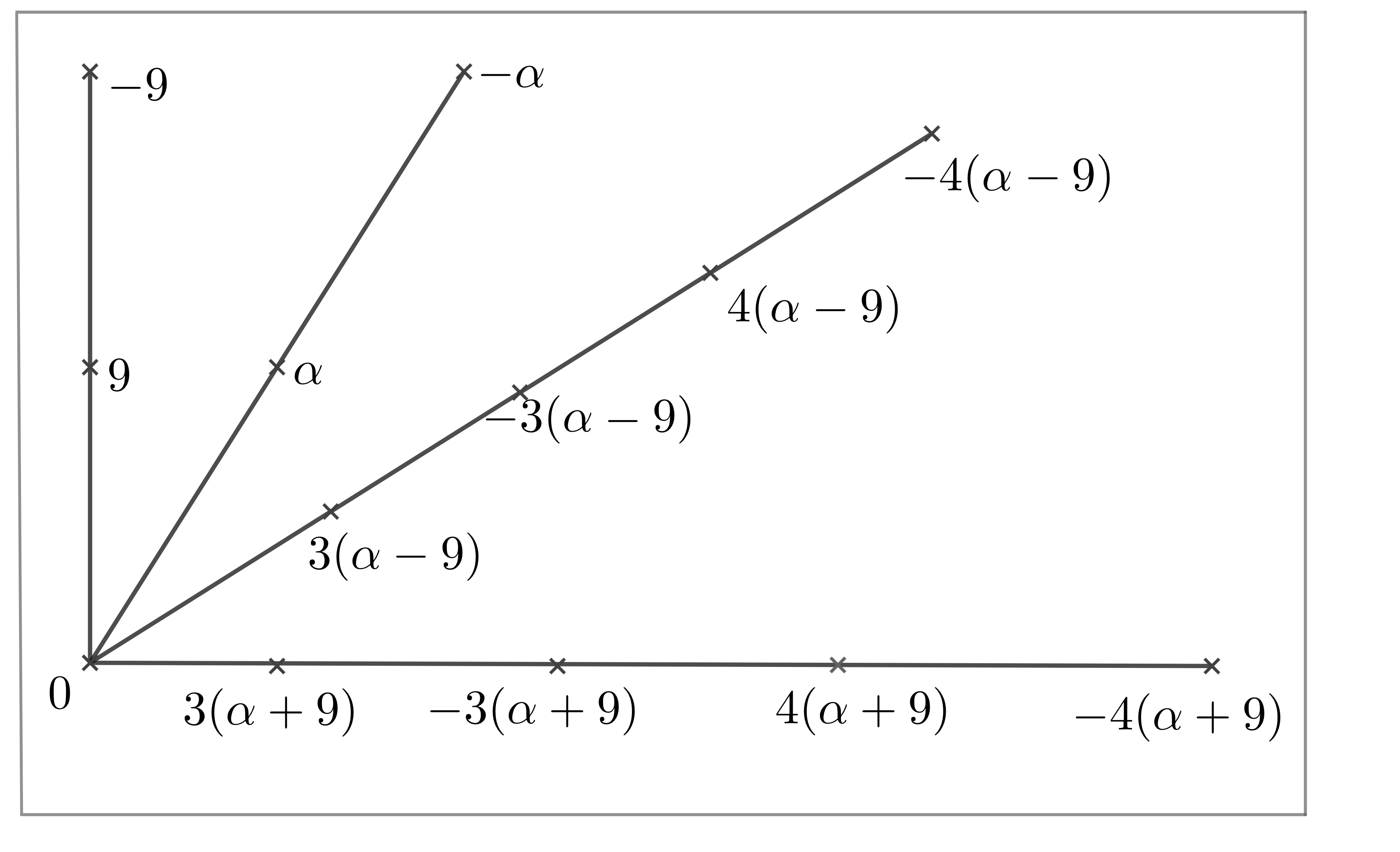}
    \quad
     \caption{\label{q=23-2}$C_{23}^{B}$} 
    \end{figure}

Let $H=\left \{ 1,-1 \right \}$ be a subset $\FF_{23}^{\ast}$, Set $C_{23}^{B}:= \{ 0 \}\cup  9H \cup \alpha H  \cup ( \alpha +9  )\{3H,4H\}\cup ( \alpha-9 )\{3H,4H\}$ be a subset of finite field $\FF_{23^2}$. By Magma, we have that $C_{23}^{B}$ is a maximal clique in $P(23^2)$ with size $\frac{q+3}{2}$ for $q=23$.
The structure of the clique $C_{23}^{B}$ is presented in Fig \ref{q=23-2}.

Remarks: There exists other lines  
that  contain at least three points in $C_{23}^{B}$. 
They are:$\{9,3(\alpha+9),4(\alpha-9)\}$, $\{9,-4(\alpha+9),-3(\alpha-9)\}$,$\{-9,-3(\alpha+9),-4(\alpha-9)\}$, $\{-9,4(\alpha+9),3(\alpha-9)\}$, $\{\alpha,3(\alpha+9),-4(\alpha-9)\}$, $\{\alpha,-4(\alpha+9),3(\alpha-9)\}$, $\{-\alpha,-3(\alpha+9),4(\alpha-9)\}$, $\{-\alpha,4(\alpha+9),-3(\alpha-9)\}$        
 
\begin{lem}
 Set $\mathcal{C}_{23}^{B}$ be the orbit of $Aut(P(23^2))$ acting on the cliques with $C_{23}^{B}\in \mathcal{C}_{23}^{B}$. Then $Aut(P(23^2))_{C_{23}^{B}}=\lg \sigma, \tau  \rg\cong \D_8$, where $ \sigma(\gamma)=18\alpha \gamma$ and $\tau(\gamma)=\gamma^{23}$ for any  $\gamma\in \FF_{23^2}$. Moreover, $|\mathcal{C}_{23}^{B}|=34914$.    
\end{lem}
\demo
Note that for any $\phi \in \Aut(P(23^2))$,
$\phi(\gamma)=a\gamma^{v}+b$, where $ a\in S$, $b\in \FF_{23^{2}}$, $v\in Gal(\FF_{23^{2}})$.  Note that the point $0$ in the clique $C_{23}^B$ is a intersecting point of two lines, which contain five points in this clique. If $\phi\in Aut(P(23^{2}))_{C_{23}^B}$, then $\phi(0)=0 $ and $b=0$. Set $\mathcal{L}:=\{9H, \alpha H\}$ be the two short lines in clique $C_{23}^{B}$,  which have a common point $0$. Then $\phi(\mathcal{L})=\mathcal{L}$. 

 If $\phi(9H)=9H$ and $\phi(\alpha H)=\alpha H$, then we have that $a\in\{1,-1\}$, $b=0$ and $v\in\{1,23\}$. 
 Set $\phi_{1}(\gamma)=\gamma$, $\phi_{2}(\gamma)=-\gamma$, $\phi_{3}(\gamma)=\gamma^{23}$ and $\phi_{4}(\gamma)=-\gamma^{23}$ for any $\gamma\in \FF_{23^2}$.
 
 If $\phi(9H)=\alpha H$ and $\phi(\alpha H)=9H$, then we have that $a\in\{18\alpha,-18\alpha\}$, $b=0$ and $v\in \{1,23\}$. Set $\phi_{5}(\gamma)=18\alpha\gamma$, $\phi_{6}(\gamma)=-18\alpha\gamma$, $\phi_{7}(\gamma)=18\alpha\gamma^{23}$ and $\phi_{8}(\gamma)=-18\alpha\gamma^{23}$ for any $\gamma\in \FF_{23^2}$. 

 Let $\sigma=\phi_{5}$, $\tau=\phi_{3}$, then we have that $\sigma^{4}=1$, $\tau^{2}=1$, $\tau^{-1}\sigma\tau=\sigma^{-1}$ then $\lg \sigma, \tau  \rg=\{ 1,\sigma,\sigma^2,\sigma^3,\tau,\tau\sigma,\tau\sigma^2,\tau\sigma^3\}=\{\phi_{1},\phi_{5},\phi_{2},\phi_{6},\phi_{3},\phi_{8},\phi_{4}, \phi_{7}\}\cong D_8$. The action of $\sigma$ and $\tau$ on the points in clique $C_{23}^{B}$ is presented in the following table.

 \begin{table}[htp]
\centering

\setlength{\tabcolsep}{0.6mm}
\begin{tabular}{cccccccccccccc}
\toprule  
$\gamma$ & 0&9&-9 & $\alpha$& $-\alpha$& $3(\alpha+9)$& $-3(\alpha+9)$  & $4(\alpha+9)$ & $-4(\alpha+9)$& $3(\alpha-9)$& $-3(\alpha-9)$  & $4(\alpha-9)$ & $-4(\alpha-9)$\\
	\midrule  
 
$\sigma(\gamma)$& 0&$\alpha$& $-\alpha$ &-9&9 & $3(\alpha-9)$& $-3(\alpha-9)$  & $4(\alpha-9)$ & $-4(\alpha-9)$&  $-3(\alpha+9)$& $3(\alpha+9)$  & $-4(\alpha+9)$ & $4(\alpha+9)$\\

$\tau(\gamma)$& 0&9&-9 & $-\alpha$& $\alpha$& $-3(\alpha-9)$& $3(\alpha-9)$  & $-4(\alpha-9)$ & $4(\alpha-9)$& $-3(\alpha+9)$& $3(\alpha+9)$  & $-4(\alpha+9)$ & $4(\alpha+9)$\\

\midrule 
\end{tabular}
\end{table}

Followed by the above arguments, we know that $\sigma(C_{23}^{B})=C_{23}^{B}$ and $\tau(C_{23}^{B})=C_{23}^{B}$, then $Aut(P(23^2))_{C_{23}^{B}}=\lg \sigma, \tau  \rg\cong \D_{8}$  where $ \sigma(\gamma)=18\alpha \gamma$ and $\tau(\gamma)=\gamma^{23}$  with $\gamma\in \FF_{23^2}$.

Note that $\left | Aut(P(23^{2})) \right |=\frac{q^2-1}{2}\times q^{2}\times 2=279312$, so we have that $|\mathcal{C}_{23}^{B}|=|Aut(P(23^2)):Aut(P(23^2))_{C_{23}^{B}}|=34914$.
\qed

\vskip 5mm
Set $C_{s,\eta,\gamma}:=\{s\eta x+\gamma| x\in C_{23}^B\cap \FF_{23^2}\}$ for $\eta\in\{1,\alpha+1,\alpha-1,\alpha+2,\alpha-2,\alpha+9\}$, where $s\in\{1\dots11\}$ and $\gamma\in\FF_{23^2}$ .

\begin{lem}

Set $\mathcal{C}_{23}^{B}$ be the orbit of $Aut(P(23^2))$ acting on the cliques with $C_{23}^{B}\in \mathcal{C}_{23}^{B}$. Then $\mathcal{C}_{23}^{B}=\bigcup_{s,\eta,\gamma}\{C_{s,\gamma,\eta } \mid s\in \{1\dots11\}, \gamma\in\FF_{23^2}, \eta\in\{1,\alpha+1,\alpha-1,\alpha+2,\alpha-2,\alpha+9\}\}$.  
\end{lem}

\demo 
It is obvious that $\bigcup_{s,\eta,\gamma}\{C_{s,\gamma,\eta } \mid s\in \{1\dots11\}, \gamma\in\FF_{23^2}, \eta\in\{1,\alpha+1,\alpha-1,\alpha+2,\alpha-2,\alpha+9\}\}\subset \mathcal{C}_{23}^B$. Now we will prove that $C_{s,\eta,\gamma}=C_{s',\eta',\gamma'}$ if and only if $s=s'$, $\gamma=\gamma'$ and $\eta=\eta'$, where $s,s'\in \{1\dots11\}$, $\gamma,\gamma'\in\FF_{{23}^2}$ and $\eta,\eta'\in\{1,\alpha+1,\alpha-1,\alpha+2,\alpha-2,\alpha+9\}$. 

Note that the intersection point of three lines in the clique $C_{s,\eta,\gamma}$ is $\gamma$. 
 If $C_{s,\eta,\gamma}=C_{s',\eta',\gamma'}$, then $\gamma=\gamma'$ and  the subset of two lines $\{9s\eta H,s\eta\alpha H\}= \{9s'\eta' H,s'\eta'\alpha H\}$, then $9s\eta\in\{9s'\eta',-9s'\eta',s'\eta'\alpha, -s'\eta'\alpha\}$. We have that $s=s'$, $\eta=\eta'$.

And then $\mathcal{C}_{23}^B=\bigcup_{\eta}\{C_{s,\gamma,\eta } \mid s\in \{1\dots11\}, \gamma\in\FF_{23^2}\}$ for $\eta\in\{1,\alpha+1,\alpha-1,\alpha+2,\alpha-2,\alpha+9\}$, because $|\bigcup_{\eta}\{C_{s,\gamma,\eta } \mid s\in \{1\dots11\}, \gamma\in\FF_{23^2}\}|=\frac{q-1}{2}\times q^{2}\times 6=34914$.
\qed

\section{Conclusion}
In this section we give a table that summarises the descriptions of the extra maximal cliques given in this paper, and formulate a problem.

Set $C$ be a representative element of the set of cliques with size $\frac{q+\epsilon}{2}$, where $\epsilon\in\{1,3\}$.
\begin{table}[htp]
\centering

\begin{tabular}{ccc}
\toprule  
$q$ & $H$ & $C$\\
	\midrule  
 
9& $ \{ 1,\delta  \}$& $ \{ 0  \} \cup H \cup   ( \alpha +\delta ^{5}  )H$\\\hline

13& $ \{ 1,3,9   \} $ &$\{ 0  \} \cup H \cup  \left ( \alpha +8 \right )H$\\\hline

13&$\{ 1,3,4 \}$ & $\{0\}\cup H\cup\{7(\alpha+1), 2(\alpha+1), 7(\alpha+7)\}$   \\\hline

17& $\left \{1,4,16,13 \right \} $ &  $\left \{ 0 \right \} \cup  H\cup \left ( \alpha +7 \right ) H$  \\\hline

17&$\{ 1,4,5 \}$ & $\{0\}\cup H\cup10(\alpha+6)H\cup\{10(\alpha+3)\}\cup 
\{6(\alpha+9)\}$  \\\hline

17& $ \{ 1,4,16,13  \} $ & $ \{ 0  \} \cup H \cup   ( \alpha +7  )\{1, 4 \}\cup    ( \alpha +10  )\{1,4  \}$  \\\hline

17& $ \{ 1,4,16,13 \} $ & $ \{ 0  \} \cup H \cup    ( \alpha +7  )\{1, 16 \}\cup  ( \alpha+10  )\{4,13  \}$ \\\hline

17& $ \{ 1,4,16,13  \} $ & $ \{0   \}\cup H \cup   ( \alpha +7  )\{ 1,4,16  \} \cup  \{ 4 ( \alpha +10 )  \}$  \\\hline

17& $\{ 1,4,16\}$  &$ \{0\}\cup H\cup (\alpha+7)\{1,4\}\cup (\alpha+10)\{1,4\}\cup \{11(\alpha+11)\}$ \\\hline

17& $\setminus$  &$ \{0\}\cup \{1,16\}\cup (\alpha+7)\{1,4\}\cup (\alpha+10)\{1,4\}\cup \{11(\alpha+11)\}\cup \{11(\alpha+6)\}$  \\\hline

11& $\{ 1,9 \}$  &$ \{0\}\cup H \cup (\alpha+5)H\cup 10(\alpha+6)H$ \\\hline

19& $\{1,-1 \}$  &$\{ 0  \}\cup   \{H,3H \} \cup 9\alpha H  \cup  ( \alpha +3 ) H \cup ( \alpha -3 ) H$  \\\hline

19& $ \{ 1,3  \}$ &$\{ 0 \}\cup  H \cup 2H ( \alpha -6 ) \cup 17H ( \alpha +6  )\cup 9H ( \alpha -4  ) \cup 10H ( \alpha +4  )$ \\\hline

23& $ \{1,3,5,17 \} $ &$ \{ 0  \} \cup H\cup 17H ( \alpha +2  ) \cup 6H ( \alpha -2  )$  \\\hline

23& $ \{ 1,-1  \}$ &$  \{ 0 \}\cup  9H \cup \alpha H  \cup ( \alpha +9  )\{3H,4H\}\cup ( \alpha-9 )\{3H,4H\}$  \\

\midrule 
\end{tabular}
\end{table}
~\\

\begin{prob}
Is there any extension of the maximal cliques described above to infinitely many values of $q \ge 25$ (these cliques would not be necessarily second largest for $q \ge 25$)?    
\end{prob}


\section*{Acknowledgments}
Huye Chen is supported by National Natural Science Foundation of China (12301446).


\begin{thebibliography}{99}








\bibitem{R.D. Baker} R.D. Baker, G.L. Ebert, J. Hemmeter, A. Woldar,
Maximal cliques in the Paley graph of square order, \emph{Journal of Statistical Planning and Inference}, {\bf56} (1996), 33-38. 

\bibitem{A. Blokhuis} A. Blokhuis,
On subsets of $GF(q^2)$ with square differences,
\emph{Indagationes Mathematicae (Proceedings)}, {\bf46} (1984), 369-372.

\bibitem{Magma} W. Bosma, J. Cannon, and C. Playoust, The Magma algebra system. I. The user language, \emph{J. Symbolic Comput.} {\bf24} (1997), 235--265.

\bibitem{S. Goryainov2} A. E. Brouwer, S. Goryainov, L. Shalaginov and C. H. Yip, Cliques in Paley graphs of square order and in Peisert graphs, in preparation.

\bibitem{P. Delsarte} P. Delsarte, An algebraic approach to the association schemes of coding theory, \emph{Philips Res. Rep.Suppl}. {\bf10} (1973).

\bibitem{Sergey Goryainov} S.V. Goryainov, V.V. Kabanov, L.V. Shalaginov, A.A. Valyuzhenich, On eigenfunctions and maximal cliques of Paley graphs of square order, \emph{Finite Fields and Their Applications}, {\bf52} (2018), 361-369.

\bibitem{S. Goryainov} S. Goryainov, A. Masley, L. Shalaginov, On a correspondence between maximal cliques in Paley graphs of square order, \emph{Discrete Mathematics},
{\bf345} (2022).

\bibitem{G.A. Jones} G.A. Jones, Paley and the Paley Graphs, in: G. Jones, I. Ponomarenko, J. {\v{S}}ir{\'a}{\v{n}}(Eds.), Isomorphisms, Symmetry and Computations in Algebraic Graph Theory, WAGT 2016, in: Springer Proceedings in Mathematics \& Statistics, vol. 305, Springer, Cham, 2020, pp. 155-183.

\bibitem{Michael Kiermaier} M. Kiermaier, S. Kurz, Maximal integral point sets in affine planes over finite fields, \emph{Discrete Mathematics}, {\bf309} (2009), 4564-4575.

\end{thebibliography}
\end{document}